\documentclass{elsarticle}
\usepackage{geometry}                
\usepackage{graphicx}
\usepackage[section]{placeins}
\usepackage{amsmath}
\usepackage{amssymb}
\usepackage{epstopdf}
\usepackage{stmaryrd}
\usepackage{pgfplots}
\usepackage{subfig}
\usepackage{tikz}
\usetikzlibrary{shapes,arrows,calc,positioning}
\usepackage[numbers]{natbib}
\usepackage{siunitx}
\usepackage{aicescover}

\pgfplotsset{
  label style={anchor=near ticklabel},
  xlabel style={yshift=0.75em},
  ylabel style={yshift=-1.em},
  tick label style={font=\tiny},
  label style={font=\scriptsize},
  legend style={font=\tiny, nodes=right},
  cycle list={{black, mark=*,mark size=1}, {red,mark=triangle*,mark size=1},{blue, mark=square*,mark size=1}, {green,mark=diamond*,mark size=1},{magenta, mark=+,mark size=1}, {black,dashed}, {black}}}

\makeatletter
\newcommand{\miniscule}{\@setfontsize\miniscule{4}{5}}
\makeatother

\newcommand{\dx}{\,{\rm d}x}

\newcommand{\ds}{\,{\rm d}\sigma}

\newcommand{\jump}[1]{\left\llbracket #1\right\rrbracket}
\newcommand{\vint}[1]{\left(#1\right)_{\mathcal{T}_h}}
\newcommand{\eint}[1]{\left<#1\right>_{\partial\mathcal{T}_h}}
\newcommand{\beint}[1]{\left<#1\right>_{\Gamma^b_h}}
\newcommand{\eints}[1]{\left<#1\right>_{\Gamma_h}}
\newcommand{\Th}{\mathcal{T}_h}
\newcommand{\XX}{\mathbb{X}}
\newcommand{\R}{\mathbb{R}}
\DeclareMathOperator{\meas}{meas}

\begin{document}

\title{Adjoint-Based Error Estimation and Mesh Adaptation for Hybridized Discontinuous Galerkin Methods}

\author[aices]{Michael Woopen\corref{cor}}
\ead{woopen@aices.rwth-aachen.de}
\cortext[cor1]{Corresponding author. Tel.: +49 241 8099135.}

\author[aices]{Georg May}
\ead{may@aices.rwth-aachen.de}

\author[igpm]{Jochen Sch\"utz}
\ead{schuetz@igpm.rwth-aachen.de}

\address[aices]{AICES Graduate School, RWTH Aachen University, Schinkelstr.\,2, 52062 Aachen, Germany}

\address[igpm]{IGPM, RWTH Aachen University, Templergraben 55, 52062 Aachen, Germany}

\aicescoverauthor{Michael Woopen, Georg May, and Jochen Sch\"utz}


\journal{Journal of Computational Physics}

\begin{abstract}
We present a robust and efficient target-based mesh adaptation methodology, building on hybridized discontinuous Galerkin schemes for (nonlinear) convection-diffusion problems, including the compressible Euler and Navier-Stokes equations. Hybridization of finite element discretizations has the main advantage, that the resulting set of algebraic equations has globally coupled degrees of freedom only on the skeleton of the computational mesh. Consequently, solving for these degrees of freedom involves the solution of a potentially much smaller system. This not only reduces storage requirements, but also allows for a faster solution with iterative solvers. The mesh adaptation is driven by an error estimate obtained via a discrete adjoint approach. Furthermore, the computed target functional can be corrected with this error estimate to obtain an even more accurate value. The aim of this paper is twofold: Firstly, to show the superiority of adjoint-based mesh adaptation over uniform and residual-based mesh refinement, and secondly to investigate the efficiency of the global error estimate.
\end{abstract}


\aicescoverpage

\maketitle


\section{Introduction}
Over the past few decades, several discretization techniques for partial differential equations have evolved. Two popular families, finite element and finite volume methods, are widely utilized in industry. However, the requirements on these techniques are continuously increasing, necessitating further work on more accurate and efficient methodologies. In view of the advantages and disadvantages of these two classes of methods, the so-called discontinuous Galerkin methods have attracted interest. Despite their popular advantages --- high-order accuracy on unstructured meshes, a variational setting, and local conservation, just to name a few --- they introduce a large number of degrees of freedom and they lead to large stencils due to an increased coupling across element interfaces (we refer to \citet{cockburn2003discontinuous} for an introduction to DG methods for various applications). Hybridization may be utilized to avoid these disadvantages. Here, the globally coupled unknowns have support on the mesh skeleton, i.e. the element interfaces, only. As a result, the global system is decreased in size and improved in terms of sparsity at the same time. The solution in the interior of the elements is then obtained using element-wise reconstruction.

Hybridization is a classic paradigm for elliptic equations~\cite{fraeijs1965displacement,arnold1985mixed}. \citet{cockburn2004characterization} gave a new characterization of the approximate solution obtained by hybridized mixed methods for second order elliptic problems. Subsequently, \citet{cockburn2009unified} introduced a unifying framework for hybridization of finite element methods for these problems. Their framework comprises hybridized mixed, continuous Galerkin, nonconforming, and hybridizable discontinuous Galerkin methods. This allows utilization of different methods in different elements, whereat the coupling is done automatically.

\citeauthor{nguyen2009implicitlin} devised a hybridizable discontinuous Galerkin method for linear and nonlinear convection-diffusion equations in \cite{nguyen2009implicitlin,nguyen2009implicitnonlin}. They observed that the approximations for the scalar variable and the flux converge with the optimal order in the $L^2$-norm. \citeauthor{peraire2010hybridizable} extended this method to the compressible Euler and Navier-Stokes equations in \cite{peraire2010hybridizable}. They showed that their method produces optimal convergence rates for both the conserved quantities and the viscous stresses and heat fluxes.


In \cite{Egger2010A-hybrid-mixed-}, \citeauthor{Egger2010A-hybrid-mixed-} proposed a hybridized discretization method for  convection-diffusion equations based on a discontinuous Galerkin discretization for convection terms, and a mixed method using H(div) spaces for the diffusive terms. This method was extended    by \citet{schutz201358} to the compressible Navier-Stokes equations. 

While the size of the global system arising from a given discretization may be reduced significantly by hybridization, it may still be considerable for many problems (for example turbulent flow around a complete airplane). Therefore, it is very important to distribute the degrees of freedom in an efficient way, resolving only those regions which are important for the values of interest. In most engineering applications, one is not necessarily interested in full flow details, but rather in some specific scalar quantities. In external aerodynamics, these may be lift or drag coefficients of airplanes. With the aim of computing accurate values for such functional quantities in the most efficient way, target-based error control methods have been developed. One such method is based on the adjoint solution of the original governing equations. In this method, an additional linear system is solved which then gives an estimate on the spatial error distribution contributing to the error in the target functional. This estimate can be used as a criterion for local mesh adaptation.
We refer to \citet{Becker2001} and \citet{giles2002adjoint} for a very comprehensive review of adjoint methods.

In the present work we follow the discrete adjoint approach.
While adjoint-based mesh adaptation is routinely applied in the context of high order methods \cite{hartmann2002adaptive,hartmann2006adaptive,fidkowski2011review}), combining hybridized methods for nonlinear advection-diffusion systems (including the compressible Euler and Navier-Stokes equations) with adjoint-based optimal control techniques has, to the best of our knowledge, not previously been attempted.

This paper is structured as follows. We will briefly cover the governing equations, namely the compressible Euler and Navier-Stokes equations, in Sec.\,\ref{sec:physics}. After that we introduce our discretization and describe the concept of hybridization in Sec.\,\ref{sec:discretization}. In Sec.\,\ref{sec:errorestimation} we establish the adjoint formulation and show how hybridization can be applied to the dual problem. Next, we validate our method with a scalar test case of which the analytical solution is known a priori and then show its efficiency and robustness with examples from compressible flow, including also the transonic regime, in Sec.\,\ref{sec:results}. At the end, we give some concluding remarks and some outlook on future work in Sec.\,\ref{sec:conclusions}.

\section{Governing Equations}
\label{sec:physics}



We consider systems of balance laws 
\begin{equation}
 \label{eq:convdiff}
 \nabla\cdot\left(f_c(w)-f_v(w,\nabla w)\right)=s\left(w,\nabla w\right)
\end{equation}
with given convective and diffusive fluxes 
\begin{equation}
 f_c:\mathbb{R}^m\rightarrow \mathbb{R}^{m\times d} ,\quad  \quad f_v:\mathbb{R}^m\times \mathbb{R}^{m\times d}\rightarrow\mathbb{R}^{m\times d}, 
\end{equation}
and a state-dependent source term 
\begin{equation}
 s:\mathbb{R}^m\times\mathbb{R}^{m\times d}\rightarrow\mathbb{R}^m. 
\end{equation}
%
We denote the spatial dimension by $d$ and the number of conservative variables by $m$.
Suitable boundary conditions are discussed below. 
For $f_v\ne 0$ or a gradient-dependent source term $s$, we may formally rewrite \eqref{eq:convdiff} as 
\begin{align}
 \label{eq:mixed_convdiff1}
 q                                       &= \nabla w \\
 \label{eq:mixed_convdiff2}
 \nabla\cdot\left(f_c(w)-f_v(w,q)\right) &=s\left(w,q\right). 
\end{align}
Formulation \eqref{eq:mixed_convdiff1}-\eqref{eq:mixed_convdiff2} is frequently applied when motivating viscous discontinuous Galerkin discretizations, as it is only a system of first-order partial differential equations. 

\subsection{Two-Dimensional Euler Equations}
The Euler equations are comprised of the inviscid compressible continuity, momentum and energy equations. 
They are given in conservative form as
\begin{equation}
\nabla\cdot f_c(w)=0,
\end{equation}
with the vector of conserved variables
\begin{equation}
w:=(\rho,\rho u, \rho v, E)^T.
\end{equation}
Here $\rho$ is the density, $u$ and $v$ are the components of the velocity vector $\widehat w:=(u,v)^T$, and $E$ is the total energy.
The convective flux is given by 
\begin{align}
f_{c,1}&=\left(\rho u, p+\rho u^2, \rho uv, u(E+p)\right)^T \\
f_{c,2}&=\left(\rho v, \rho uv, p+\rho v^2, v(E+p)\right)^T.
\end{align}

Pressure is related to the conservative flow variables $w$ by the equation of state
\begin{equation}
p=(\gamma-1)\left(E-\frac12\rho\left(u^2+v^2\right)\right)
\end{equation}
where $\gamma=c_p/c_v$ is the ratio of specific heats, generally taken as $1.4$ for air. 

Along wall boundaries we apply the slip boundary condition
\begin{equation}
\widehat{w}_n(w) := (u,v)\cdot n=0.
\end{equation}
For later use, we define a boundary function $w_{\partial \Omega}(w)$ which satisfies $\widehat{w}_n(w_{\partial\Omega}(w))=0$. Such a function can  be written
\begin{equation}
\label{eq:boundary_euler}
w_{\partial\Omega}(w):=\begin{pmatrix}
1 & 0 & 0 & 0\\
0 & 1-n^2_x & -n_xn_y & 0\\
0 & -n_xn_y & 1-n^2_y & 0\\
0 & 0 & 0 & 1
\end{pmatrix}w
\end{equation}
where $n=(n_x,n_y)$.

\subsection{Two-Dimensional Navier-Stokes Equations}
The Navier-Stokes equations in conservative form  are given by
\begin{equation}
\nabla\cdot\left(f_c(w)-f_v(w,\nabla w)\right)=0.
\end{equation}
The convective part $f_c$ of the Navier-Stokes equations coincides with the one of the Euler equations. The viscous flux is given by 
\begin{align}
f_{v,1}&=\left(0,\tau_{11},\tau_{21},\tau_{11}u+\tau_{12}v+kT_x\right)^T \\
f_{v,2}&=\left(0,\tau_{12},\tau_{22},\tau_{21}u+\tau_{22}v+kT_y\right)^T.
\end{align}
The temperature is defined via the ideal gas law
\begin{equation}
T=\frac{\mu\gamma}{k\cdot\mathrm{Pr}}\left(\frac{E}{\rho}-\frac12\left(u^2+v^2\right)\right)=\frac{1}{(\gamma-1)c_v}\frac{p}{\rho}
\end{equation}
where $\mathrm{Pr}=\frac{\mu c_p}{k}$ is the Prandtl number, which for air at moderate conditions can be taken as a constant with a value of $\mathrm{Pr}\approx 0.72$. $k$ denotes the thermal conductivity coefficient.
For a Newtonian fluid, the stress tensor is defined as 
\begin{equation}
\tau=\mu\left(\nabla\widehat{w}+\left(\nabla\widehat{w}\right)^T-\frac23\left(\nabla\cdot\widehat{w}\right)\mathrm{Id}\right).
\end{equation}

The variation of the molecular viscosity $\mu$ as a function of temperature is determined by Sutherland's law as
\begin{equation}
\mu=\frac{C_1T^{3/2}}{T+C_2}
\end{equation}
with $C_1=\SI{1.458e-6}{\frac{kg}{ms\sqrt{K}}}$ and $C_2=\SI{110.4}{K}$.

Along wall boundaries, we apply the no-slip boundary condition, i.e.
\begin{equation}
\widehat{w}(w):=(u,v)=0
\end{equation}
with corresponding boundary function
\begin{equation}
\label{eq:boundary_ns}
w_{\partial\Omega}(w):=\left(\rho, 0, 0, E\right)^T
\end{equation}
so that $\widehat{w}(w_{\partial\Omega}(w))=0$.
Furthermore, one has to give boundary conditions for the temperature. In the present work we use the adiabatic wall condition, i.e.
\begin{equation}
\nabla T\cdot n=0.
\end{equation}
Combining both no-slip and adiabatic wall boundary conditions, gives a condition for the viscous flux, namely
\begin{equation}
\label{eq:boundary_flux_ns}
f_{v,\partial\Omega}(f_v)=\begin{pmatrix}
0 & \tau_{11} & \tau_{21} & 0 \\
0 & \tau_{12} & \tau_{22} & 0
\end{pmatrix}^T.
\end{equation}


\section{Discretization}
\label{sec:discretization}

\subsection{Notation}

In order to discretize equation \eqref{eq:convdiff}, we consider a triangulation  $\Th = \{ K \}$ of an open set $\Omega\subset \mathbb{R}^2$, such that 
\begin{align}
\overline \Omega = \bigcup_{K \in \Th } K,
\end{align}
and the intersection of $K, K' \in \Th$ is either the empty set or a lower-dimensional manifold. 
In this publication, we use only triangular elements. To define the HDG method, we need two definitions of \emph{edge boundaries}, the distinction will become clearer in the sequel. We define the set of element boundaries $\partial \Th$ and the skeleton of the mesh $\Gamma_h$ to be 
\begin{alignat}{3}
 \partial \Th &:= \{ \ \partial K \backslash \partial \Omega \ :\  K \in \Th \ \}, &\quad& \\
 \Gamma_h     &:= \{ \ e \ : \ e = K \cap K'  \text{ for } K, K' \in \Th; \meas_{d-1}(e)\neq 0 \ \}.
\end{alignat}
Note that elements in $\Gamma_h$ appear twice in $\partial \Th$ with different orientation. Note furthermore that both $\Gamma_h$ and $\partial \Th$ do not include boundary edges. The set of boundary edges is denoted by $\Gamma_h^b$. Based on these sets, we can introduce both inner and edge-oriented scalar products by 
\begin{align*}
  \vint{v,w}&:=\sum_{K\in\mathcal{T}_h}\int_K vw\dx, \\ 
  \eint{v,w}&:=\sum_{K\in\mathcal{T}_h}\int_{\partial K}vw\ds,\qquad
  \eints{v,w}:=\sum_{e\in\Gamma_h}\int_e\ vw\ds.
\end{align*}
The scalar product $\left<\cdot, \cdot \right>_{\Gamma_h^b}$ is defined analogously. 
The algorithm to be presented will approximate both the unknown $w$ and its gradient $q$ (see \eqref{eq:mixed_convdiff1}-\eqref{eq:mixed_convdiff2}), as well as the projection of $w$ onto the skeleton of the mesh, i.e., the function $\lambda := w_{|\Gamma_h}$. To this end, we introduce the function spaces
%
\begin{alignat}{5}
V_h &:= \{v\in L^2\left(\Omega\right)     &\ :\ & v\vphantom{\big|}_{|_K}     &\ \in\ & \Pi^{p}(K), &\quad& K&\ \in\ &\mathcal{T}_h\}^{m\times d}\\
W_h &:= \{w\in L^2\left(\Omega\right)     &\ :\ & w\vphantom{\big|}_{|_K}     &\ \in\ & \Pi^{p}(K), &\quad&K&\ \in\ &\mathcal{T}_h\}^m\\
M_h &:= \{\mu\in L^2\left(\Gamma_h\right) &\ :\ & \mu\vphantom{\big|}_{|_e} &\ \in\ & \Pi^{p}(e), &\quad&e&\ \in\ &\Gamma_h\}^m,
\end{alignat}
where $ \Pi^{p}$ is the space of polynomials of degree at most $p$. 
The hybridized discontinuous Galerkin method seeks approximations $(q_h, w_h, \lambda_h) \in V_h \times W_h \times M_h$. For convenience, we introduce the abbreviation 
\begin{align}
 \XX_h :=V_h \times W_h \times M_h. 
\end{align}
Furthermore, we define the number of elements and (interior) faces, respectively, by 
\begin{align}
 n_e := |\Th|, \quad n_f := |\Gamma_h|.
\end{align}
\subsection{Method}

A classical discontinuous Galerkin discretization starts from the strong form of \eqref{eq:mixed_convdiff1}-\eqref{eq:mixed_convdiff2}, tests it against $\left(\tau_h,\varphi_h\right)\in V_h\times W_h$ and introduces numerical fluxes $\widehat w$, $\widehat {f}_c$ and $\widehat f_v$ to obtain 
\begin{align}
\label{eq:dg1}
0=&\phantom{-}\vint{\tau_h,q_h}+\vint{\nabla\cdot\tau_h,w_h}-\eint{\tau_h\cdot n,\widehat{w}}\\
\label{eq:dg2}
&-\vint{\nabla\varphi_h,f_c(w_h)-f_v(w_h,q_h)}-\vint{\varphi_h, s(w_h,q_h)}+\eint{\varphi_h, \widehat{f_c}-\widehat{f_v}}\\
\label{eq:dg3}
&+\mathcal{N}^{\rm DG}_{h,\partial\Omega}\left(q_h,w_h;\tau_h,\varphi_h\right)+\mathcal{N}_{h,\rm sc}\left(q_h,w_h;\varphi_h\right). 
\end{align}
In order to realize a consistent and stable method, a suitable choice  of the fluxes is as important as the appropriate choice of the shock-capturing operator $\mathcal{N}_{h,\rm sc}\left(q_h,w_h;\varphi_h\right)$ and the boundary operator $\mathcal{N}^{\rm DG}_{h,\partial\Omega}\left(q_h,w_h;\tau_h,\varphi_h\right)$, see, e.g., \cite{UnifiedAnalysisDG}. In the context of DG methods, $\widehat w \in M_h$ is a function of $w^-$ and $w^+$. It has been realized, already in the context of mixed methods \cite{arnold1985mixed}, that one can approximate this quantity independently of $w$. To ensure unique solvability of the whole system, one has to introduce an additional equation, which guarantees that the total flux $\widehat f_c + \widehat f_v$ has a (weak) divergence, i.e., is continuous over the edges in normal direction. 
This insight lies at the core of the idea of the HDG method. More precisely, introducing $\widehat w := \lambda_h$ into \eqref{eq:dg1}-\eqref{eq:dg3} and adding the continuity equation for $\widehat f_c + \widehat f_v$, the HDG method seeks solutions $(q_h,w_h,\lambda_h)\in \XX_h$ s.t. $\forall(\tau_h,\varphi_h,\mu_h)\in \XX_h$
\begin{align}
\label{eq:hdg1}
0=&\phantom{-}\mathcal{N}_h\left(q_h,w_h,\lambda_h;\tau_h,\varphi_h,\mu_h\right)\\
\label{eq:hdg2}
:=&\phantom{-}\vint{\tau_h,q_h}+\vint{\nabla\cdot\tau_h,w_h}-\eint{\tau_h\cdot n,\lambda_h} \\
\label{eq:hdg3}
&-\vint{\nabla\varphi_h,f_c(w_h)-f_v(w_h,q_h)}-\vint{\varphi_h, s(w_h,q_h)}+\eint{\varphi_h, \widehat{f_c}-\widehat{f_v}}  \\
\label{eq:hdg3b}
&+\mathcal{N}_{h,\partial\Omega}\left(q_h,w_h;\tau_h,\varphi_h\right)+\mathcal{N}_{h,\rm sc}\left(q_h,w_h;\varphi_h\right)\\
\label{eq:hdg4}
&+\eints{\mu_h,\jump{\widehat{f_c}-\widehat{f_v}}} .
\end{align}

The choice of fluxes is motivated by a Lax-Friedrich's numerical flux \cite{laxfriedrichs} for the convective part, and a local discontinuous Galerkin flux \cite{localdg} for the viscous part. More precisely, we choose $\widehat f_c$ and $\widehat f_v$ by 
\begin{alignat}{2}
\widehat{f_c}\left(\lambda_h,w_h\right)&=f_c\left(\lambda_h\right)\cdot n          &\ -\ &\alpha_c\left(\lambda_h-w_h\right) \\
\widehat{f_v}\left(\lambda_h,w_h, q_h\right)&=f_v\left(\lambda_h, q_h\right)\cdot n&\ +\ &\alpha_v\left(\lambda_h-w_h\right),  
\end{alignat}
where both $\alpha_c$ and $\alpha_v$ are (strictly) positive real-valued parameters. Note that they could equally well be given by a tensor, however, we rely on scalars. Furthermore, it is worth mentioning that the choice of these parameters is not very sensitive to the quality of the resulting algorithm. 

Note that the use of $\eint{\cdot,\cdot}$ in \eqref{eq:hdg2}--\eqref{eq:hdg3b} and the use of $\eints{\cdot, \cdot}$ in \eqref{eq:hdg4} has a direct application in the implementation of the HDG method: \eqref{eq:hdg2}--\eqref{eq:hdg3} can be solved locally for $w_h$ and $\sigma_h$, and are thus called \emph{local solvers}. \eqref{eq:hdg4} is solved globally, but only in terms of $\lambda_h$. 

It has been recognized that in order to use the \emph{discrete adjoint} approach to discretize the continuous adjoint equations, one has to use an \emph{adjoint consistent} scheme \cite{OliverDarmofalAdjointConsistency, hartmann2007adjoint}. The crucial part here is the discretization of boundary terms. In \cite{schutz2012adjoint}, Sch\"utz and May have investigated the concept of adjoint consistency for general hybrid mixed methods and showed how to discretize the boundary terms. We adopt this methodology and discretize the boundary terms using boundary states $w_{\partial\Omega}$ and boundary fluxes $f_{v,\partial\Omega}$ imposed on $w$ and $q$. (For a definition of these quantities, please refer to \eqref{eq:boundary_euler}, \eqref{eq:boundary_ns} and \eqref{eq:boundary_flux_ns}.) More precisely, we choose 
\begin{align*}
\mathcal{N}_{h,\partial\Omega}\left(q_h,w_h;\tau_h,\varphi_h\right):=&\phantom{+}\beint{\tau_h\cdot n,w_{\partial\Omega}\left(w_h\right)}\\
&+\beint{\varphi_h,\left(f_c\left(w_{\partial\Omega}\left(w_h\right)\right)-f_{v,\partial\Omega}\left(f_v\left(w_{\partial\Omega}\left(w_h\right),q_h\right)\right)\right)\cdot n}.
\end{align*}
Note that in contrast to the methods considered in \cite{schutz2012adjoint}, we do not use $\lambda_h$ here. (In fact, as $\lambda_h \in M_h$, $\lambda_h$ is not defined at the boundary.) However, the analysis can easily be extended to show that the resulting method is adjoint consistent. 

We define the shock-capturing terms in a standard manner~\cite{jaffre:1995convDG}
\begin{equation}
\mathcal{N}_{h,\rm sc}\left(q_h,w_h;\varphi_h\right):=\vint{\nabla\varphi_h,\epsilon\left(w_h,q_h\right)q_h}, 
\end{equation}
however using a diffusion coefficient $\epsilon\left(w, \nabla w\right)$ which avoids coupling between elements, in order to maintain the locality of the local solves. We use the sensor developed by Nguyen and Peraire \cite{nguyen2011adaptive} that is based on the local dilatation of the flow.
In case we  approximate the Euler equations, we replace $q_h$ by $\nabla w_h$.
%


%

\subsection{Relaxation}
To simplify the notation, we define the vectors 
\begin{align}
 {\bf x}_h:=\left(q_h,w_h,\lambda_h\right) , \quad {\bf y}_h:=\left(\tau_h,\varphi_h,\mu_h\right) 
\end{align}
and abbreviate \eqref{eq:hdg1} by 
\begin{align}
 \label{eq:hdg_short}
 \mathcal{N}_h\left({\bf x}_h, {\bf y}_h\right) &= 0 \quad \forall{\bf y}_h \in \XX_h.
\end{align}
The HDG method \eqref{eq:hdg_short} constitutes a nonlinear system of equations that we solve by means of a damped Newton-Raphson method. To this end, starting from an initial guess ${\bf x}^0_h := \left(q_h^0,w_h^0,\lambda_h^0\right)$, we iteratively solve the linearized system 
\begin{equation}
\vint{\varphi_h,\frac1{\Delta t_K^n}\delta{w}^n_h} + \mathcal{N}^{\prime}_h\left[{\bf x}^n_h\right]\left({\delta{\bf x}}^n_h;{\bf y}_h\right)=-\mathcal{N}_h\left({\bf x}^n_h;{\bf y}_h\right) \quad \forall {\bf y}_h \in \XX_h
\label{eq:linearsystem}
\end{equation}
for ${\delta{\bf x}}^n_h \in \XX_h$. 
The solution is updated in the usual way as 
\begin{equation}
{\bf x}^{n+1}_h={\bf x}^n_h+{\delta{\bf x}}^n_h.
\end{equation}
This update step is repeated until the residual $\mathcal{N}_h\left({\bf x}^n_h;{\bf y}_h\right)$ drops below a prescribed threshold. 
$\mathcal{N}^{\prime}_h$ denotes the Fr\'echet derivative of $\mathcal{N}_h$ with respect to the first argument. 

The damping term that resembles a time-derivative as it occurs in the implicit Euler method is introduced to increase the radius of convergence of the iteration. Usually, the damping parameter $\Delta t_K^n$ is small at the beginnning of the iteration, and tends to infinity in dependence of the quality of ${\bf x}_h^{n}$ to obtain the full Newton-Raphson method. No time accuracy is needed, so we choose this parameter cell-dependent as 
\begin{equation}
\Delta t^n_K:={\rm CFL}^n\frac{\left|K\right|}{\lambda_c+4\lambda_v}, 
\end{equation}
with $\lambda_c$ and $\lambda_v$ representing approximations to the maximum eigenvalue of the Jacobians of convective and diffusive flux, respectively (see \citet{mavriplis1990multigrid} for a more detailed presentation). ${\rm CFL}^n$ is a prescribed scalar that depends on ${\bf x}^{n}_h$ or rather on the residual.

\subsection{Hybridization}\label{sec:hybrid_primal}
Equation \eqref{eq:linearsystem} is a linear equation in $\delta {\bf x}_h^n$. As such, it constitutes a large system of equations. However, the fact that $\delta q_h^n$ and $\delta w_h^n$ are only implicitly (via $\delta \lambda_h^n$) coupled over element boundaries is exploited in this section to reduce the dimension of the system. 

We assume that $\left[\delta Q, \delta W, \delta\Lambda\right]$ are basis coefficients corresponding to $\delta {\bf x}_h^n$. Then, obviously, \eqref{eq:linearsystem} is equivalent to 
\begin{equation}
\left[\begin{array}{ccc}A & B & R \\C & D & S \\ L & M & N\end{array}\right]\left[\begin{array}{c}\delta Q \\\delta W \\ \delta\Lambda\end{array}\right]=\left[\begin{array}{c}F \\ G \\ H\end{array}\right],
\end{equation}
which can be split into 
\begin{equation}
\label{eq:ls}
\left[\begin{array}{cc}A & B \\C & D\end{array}\right]\left[\begin{array}{c}\delta Q \\\delta W\end{array}\right]=\left[\begin{array}{c}F \\ G\end{array}\right]-\left[\begin{array}{c}R \\ S\end{array}\right]\delta\Lambda
\end{equation}
and
\begin{equation}
\label{eq:cons}
\left[\begin{array}{cc}L & M\end{array}\right]\left[\begin{array}{c}\delta Q \\\delta W\end{array}\right]+N\delta\Lambda=H.
\end{equation}
Eq.~\eqref{eq:ls} is the part that corresponds to the local solution processes \eqref{eq:hdg2}-\eqref{eq:hdg3}. As such, matrices $A$, $B$, $C$ and $D$ are block-diagonal and can be easily inverted. Therefore, we can explicitly solve \eqref{eq:ls} and insert the expressions for $\left[\delta Q, \delta W\right]$ into \eqref{eq:cons} to obtain a system in $\delta \lambda_h$ only. It is straightforward to show that this system is given by 
\begin{equation}
\label{eq:hybridsystem}
\left(N-\left[\begin{array}{cc}L & M\end{array}\right]\left[\begin{array}{cc}A & B \\C & D\end{array}\right]^{-1}\left[\begin{array}{c}R \\S\end{array}\right]\right)\delta\Lambda=H-\left[\begin{array}{cc}L & M \end{array}\right]\left[\begin{array}{cc}A & B \\C & D\end{array}\right]^{-1}\left[\begin{array}{c}F \\ G\end{array}\right].
\end{equation}
The matrix on the left-hand side of \eqref{eq:hybridsystem} is a sparse $ n_f\times n_f$ block matrix. The number of values in each block is of the order $\mathcal O\left( (m p^{d-1})^2 \right)$. Any block row consists of one dense diagonal block and (for simplices) $2d$ dense off-diagonal blocks, so in total, the number of nonzero entries is of the order $\mathcal O \left( (2d+1) n_f m^2 p^{2(d-1) } \right)$. (Note that an edge $e = K \cap K'$ of a triangle has four neighboring edges, as the local solution process affects all the edges of both $K$ and $K'$.) In comparison, the matrix corresponding to a DG discretization is also a sparse block matrix, its number of nonzeros entries is, however, of the order $\mathcal O \left((d + 2) n_e m^2 p^{2d} \right)$. For large $p$, this usually requires more nonzero values to be stored. This constitutes an advantage of an HDG method both for memory requirements and less floating point operations. 

The system \eqref{eq:hybridsystem} is solved with an ILU(n) preconditioned GMRES available through the PETSc library \cite{petsc-web-page, petsc-user-ref}. Furthermore, we want to note that we use the Netgen/Ngsolve library \cite{schoberl1997netgen} which offers besides geometry handling and mesh generation also quadrature rules and the evaluation of basis functions for numerous finite elements.


\section{Error Estimation and Adaptation}
\label{sec:errorestimation}

\subsection{Adjoint-Based Error Estimation}
In engineering applications, it is usually only a few real-valued quantity that are of interest, and not necessarily the solution quality per se. In aerodynamic applications, these quantities can e.g. be lift and drag, in two dimensions given by 
\begin{align}
 \label{eq:liftdrag}
 J(q, w) &:= \left<p(w) \beta - \tau(q, w) \beta, n \right>_{\partial\Omega} 
\end{align}
for suitable choices of $\beta \in \R^2$. Note that integration is only along the airfoil. 
In the HDG context, the target functional $J$ is discretized by $J_h$ which is evaluated with exact boundary operators $w_{\partial \Omega}$ and $f_{v, \partial \Omega}$.  For details, we refer to \cite{schutz2012adjoint}. 

Adjoint-based error estimation has been developed in order to most accurately approximate such quantities \cite{giles2002adjoint,Becker2001}. To put this into a mathematical context, let a functional $J_h:\XX_h\rightarrow\mathbb{R}$ be given and consider the error in this quantity, 
\begin{equation}
e_h:=J_h\left({\bf x}\right)-J_h\left({\bf x}_h\right). 
\end{equation}
Formally, the derivation of the (discrete) adjoint equation relies on Taylor's expansion of both $J_h$ and $\mathcal{N}_h$: 
\begin{align}
\label{eq:linerr}
J_h\left({\bf x}\right)-J_h\left({\bf x}_h\right) &=J_h^ {\prime}\left[{\bf x}_h\right]\left({\bf x}-{\bf x}_h\right)+\mathcal{O}\left(\|{\bf x}-{\bf x}_h\|^2\right) \\ 
\label{eq:linres}
\mathcal{N}_h\left({\bf x}; {\bf y}_h\right)-\mathcal{N}_h\left({\bf x}_h;{\bf y}_h\right)&=\mathcal{N}_h^ {\prime}\left[{\bf x}_h\right]\left({\bf x}-{\bf x}_h;{\bf y}_h\right)+\mathcal{O}\left(\|{\bf x}-{\bf x}_h\|^2\right).
\end{align}
Obviously, the contribution $\mathcal{N}_h\left({\bf x}; {\bf y}_h\right)$ is zero, as we are using a consistent method. 
Now upon assuming that we can solve the discrete adjoint equations 
\begin{align}
 \label{eq:discreteadjoint}
\mathcal{N}^{\prime}_h\left[{\bf x}_h\right]\left({\bf y};{\bf z}\right)=J_h^{\prime}\left[{\bf x}_h\right]\left({\bf y}\right)\qquad\forall{\bf y}\in{\XX},
\end{align}
we obtain, combining \eqref{eq:linerr} and \eqref{eq:linres} (with ${\bf x}, {\bf x}_h \in \XX$, i.e., $\XX_h \subset \XX$)
\begin{equation}
e_h = -\mathcal{N}_h\left({\bf x}_h;{\bf z}\right) +\mathcal{O}\left(\|{\bf x}-{\bf x}_h\|^2\right).
\label{eq:errest_exact}
\end{equation}
As such, the adjoint solution ${\bf z}$ constitutes the connection between variations in the residual and variations in the target functional. 
Of course, besides being linear, \eqref{eq:discreteadjoint} is as difficult to solve as the original convection-diffusion equation  \eqref{eq:convdiff}, so we approximate it by seeking ${\bf z}_h \in \widetilde \XX_h$ such that 
\begin{equation}
\mathcal{N}^{\prime}_h\left[{\bf x}_h\right]\left({\bf y}_h;{\bf z}_h\right)=J_h^{\prime}\left[{\bf x}_h\right]\left({\bf y}_h\right)\qquad\forall{\bf y}_h\in\widetilde{\XX}_h,
\label{eq:adjoint}
\end{equation}
for some function space $\widetilde \XX_h \supset \XX_h$. An error estimator is obtained upon substituting ${\bf z}_h$ into \eqref{eq:errest_exact} to obtain 
\begin{equation}
e_h \approx \eta := -\mathcal{N}_h\left({\bf x}_h;{\bf z}\right).
\label{eq:errest}
\end{equation}
Note that the choice $\widetilde \XX_h = \XX_h$ would yield, due to Galerkin orthogonality, a useless error estimator. Therefore, $\widetilde \XX_h$ is chosen as a superset of $\XX_h$. It can, e.g., be obtained via mesh-refinement or increase of polynomial degree. In the DG and HDG context, the latter is more advantageous with respect to both implementation and efficiency. 

\eqref{eq:errest} can be applied to improve the quality of the target output $J_h({\bf x}_h)$, as 
\begin{align}
\label{eq:improved_J}
J_h({\bf x}) \approx J_h({\bf x}_h) + \eta.
\end{align}
 Furthermore, in its localized version 
\begin{equation}
\label{eq:local_eta}
\eta_K:=\left|\mathcal{N}_h\left({\bf x}_h;{\bf z}_h\vphantom{\big|}_{|_K}\right)\right|
\end{equation}
it can be used to drive an adaptation process. 

\subsection{Hybridization}
It is well-known, see, e.g., \cite{hartmann2007adjoint} that, in order to compute the discrete adjoint of a DG discretization, one can take the \emph{transpose} of the Jacobian matrix evaluated at $w_h$ and solve the corresponding linear system of equations with a right-hand side that depends on the target functional under consideration. This makes it very easy to implement these methods, as the main implementational effort goes into computing the Jacobian. In this section, we show that one has a similar structure in the context of HDG methods. 

We assume that $\left[\widetilde{Q}, \widetilde{W}, \widetilde{\Lambda}\right]$ are the basis coefficients corresponding to ${\bf z}_h=\left(\widetilde{q}_h, \widetilde{w}_h, \widetilde{\lambda}_h\right)\in\widetilde{\XX}_h$. Consequently, the adjoint equations \eqref{eq:adjoint} can be written as
\begin{equation}
\left[\begin{array}{ccc}A & B & R \\C & D & S \\ L & M & N\end{array}\right]^T\left[\begin{array}{c}\widetilde{Q} \\\widetilde{W} \\ \widetilde{\Lambda}\end{array}\right]=\left[\begin{array}{c}\widetilde{F} \\ \widetilde{G} \\ 0\end{array}\right].
\end{equation}
Note that there is a zero vector in the right-hand side, as there is no contribution from $\lambda_h$ to the target functional (see for example Eq.~\eqref{eq:liftdrag}). In the same spirit as in Sec. \ref{sec:hybrid_primal}, we can perform a static condensation to obtain a linear system in $\widetilde \Lambda$ viz
\begin{equation}
\label{eq:adjointhybridsystem}
\small
\left(N-\left[\begin{array}{cc}L & M\end{array}\right]\left[\begin{array}{cc}A & B \\C & D\end{array}\right]^{-1}\left[\begin{array}{c}R \\S\end{array}\right]\right)^T\widetilde{\Lambda}=-\left[\begin{array}{cc}R^T & S^T\end{array}\right]\left[\begin{array}{cc}A & B \\C & D\end{array}\right]^{-T}\left[\begin{array}{c}\widetilde{F} \\ \widetilde{G}\end{array}\right]. 
\end{equation}
Comparing \eqref{eq:adjointhybridsystem} to \eqref{eq:hybridsystem}, one can see that also in this context, the adjoint matrix is the transpose of the forward Jacobian. The same holds true for computing the local adjoint problems by 
\begin{equation}
\label{eq:adjoint_local}
\left[\begin{array}{cc}A & B \\C & D\end{array}\right]^T\left[\begin{array}{c}\widetilde{Q} \\\widetilde{W}\end{array}\right]=\left[\begin{array}{c}\widetilde{F} \\ \widetilde{G}\end{array}\right]-\left[\begin{array}{cc}L & M\end{array}\right]^T\widetilde{\Lambda}, 
\end{equation}
compare to \eqref{eq:ls}.

\subsection{Element Marking and Summary}
\label{sec:marking}
The localized error estimator \eqref{eq:local_eta} is used to drive the adaptation process. We mark elements for refinement via so-called fixed-fraction approach \cite{fidkowski2011review}, where a user-defined fraction $\theta \in (0, 1)$ is chosen. 
Arranging the elements $K_i, i = 1, \ldots, N$ in such a way that $\eta_{K_i}$ is non-increasing, we refine the elements $K_i, i = 1, \ldots ,\lfloor\theta N\rfloor$, thus obtaining a new triangulation $\mathcal{T}_{h'}$. 

\begin{figure}[htb]
    \begin{center}
    \includegraphics{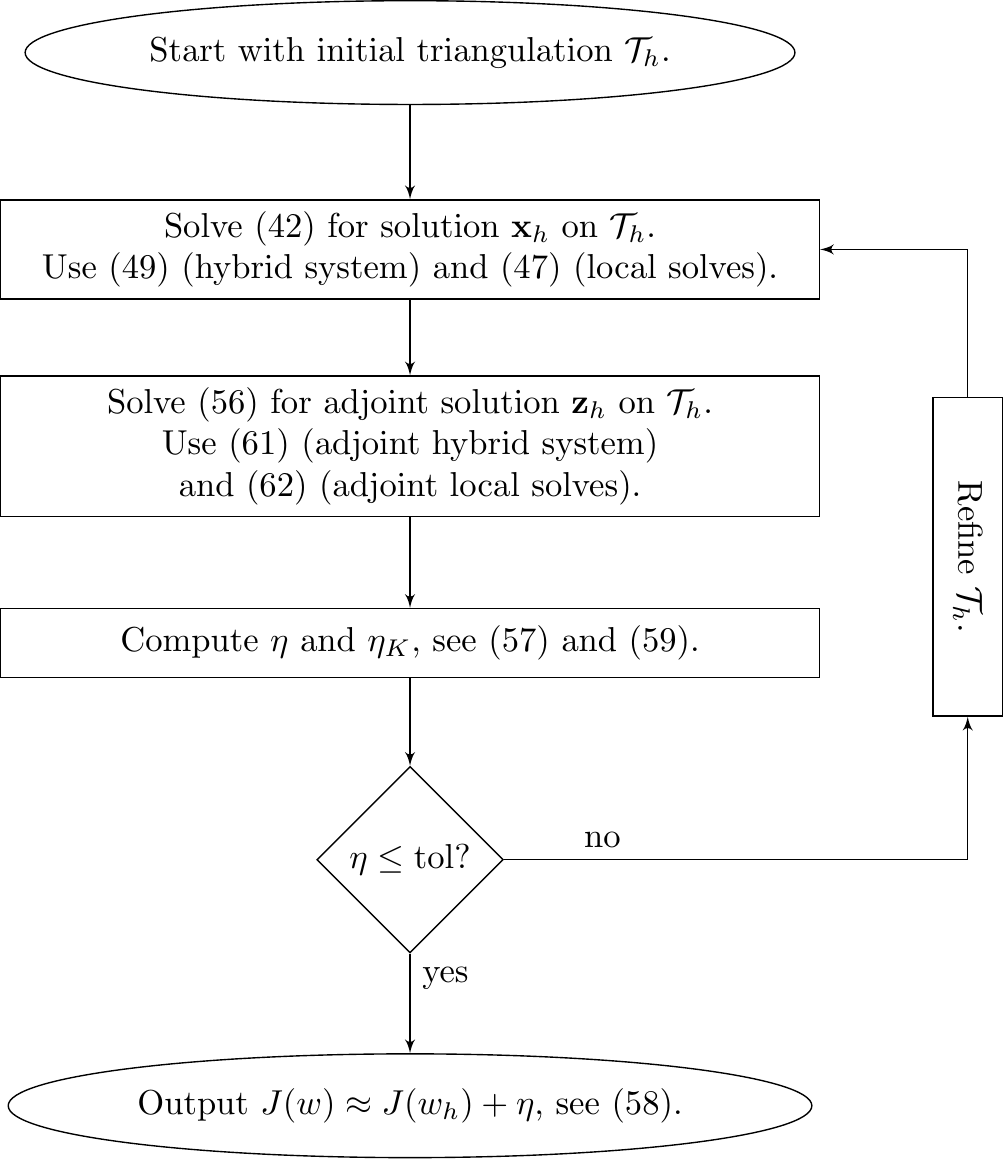}
    \caption{Summary of the total adaptive algorithm. $\text{tol}$ is a user-defined tolerance. }
    \label{fig:flowchart}
    \end{center}
\end{figure}

A summary of the overall algorithm can be found in Fig. \ref{fig:flowchart}. 

%


\section{Numerical Results}
\label{sec:results}
We validate our method with a scalar example and then move on to problems of compressible fluid flow. 
The latter include inviscid flow over a smooth bump and both inviscid and viscous flow around the NACA~0012 airfoil. The specific flow conditions are taken from the high-order workshop (see \citet{FLD:FLD3767}).

For each test case we perform adjoint- and residual-based refinement where we found $\theta\approx 0.1\textrm{--}0.2$ to be a robust choice. Furthermore, we conduct a uniform refinement study for the first two cases. For admissible target functionals and smooth primal and dual solution we expect convergence orders of $2p+1$ for both inviscid and viscous computations. This is due to the mixed formulation and in contrast with normal DG discretizations where the gradient exhibits suboptimal convergence so that the expected order for viscous target functionals is only $2p$.

\subsection{Convection-Diffusion}
\begin{figure}[h]
\centering
  \subfloat[Primal solution ($p=2$)\label{fig:bdrylayer_primal}]{\includegraphics[width=0.4\textwidth]{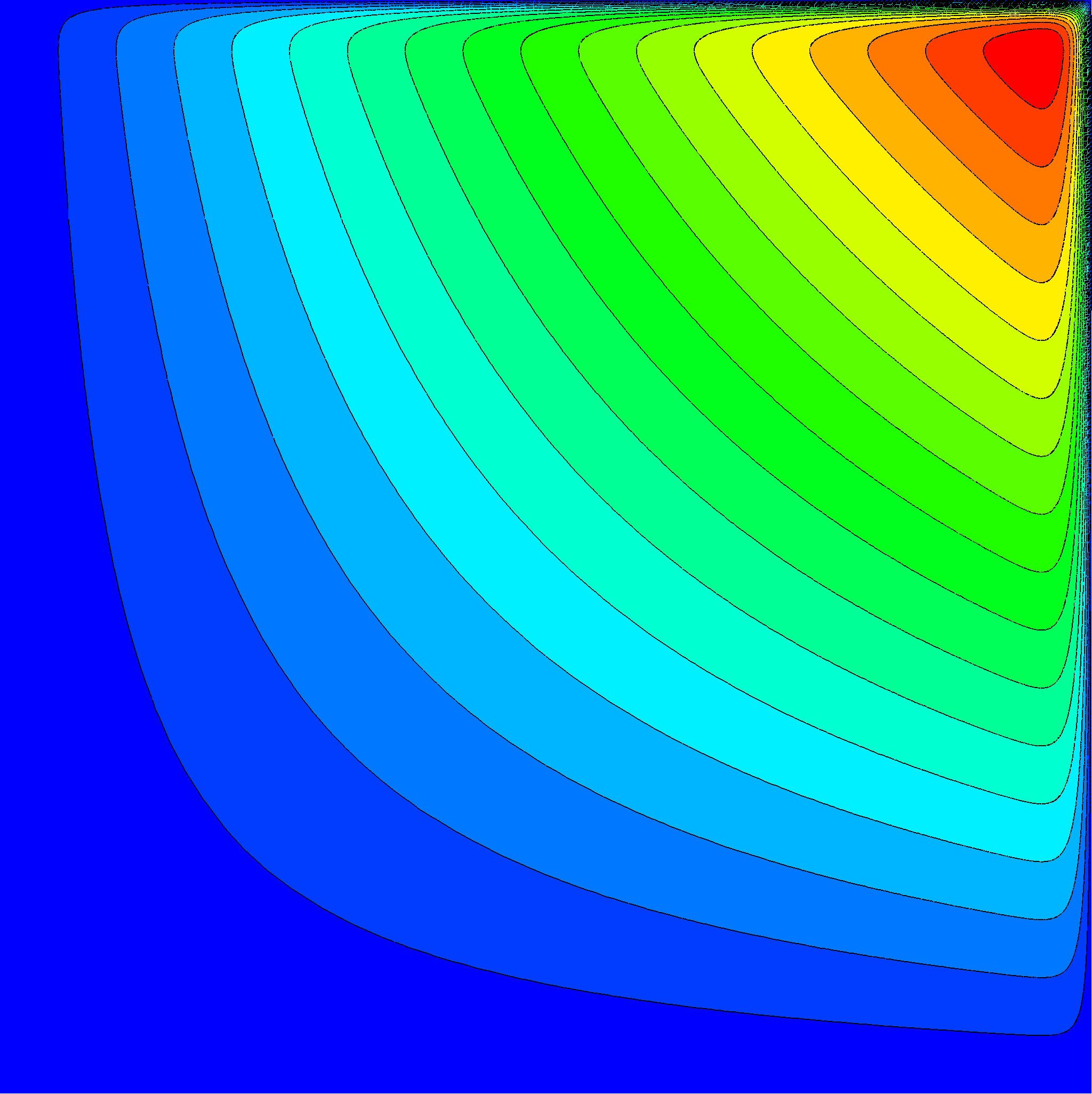}}
  \qquad
  \subfloat[Dual Solution ($p=3$)\label{fig:bdrylayer_dual}]{\includegraphics[width=0.4\textwidth]{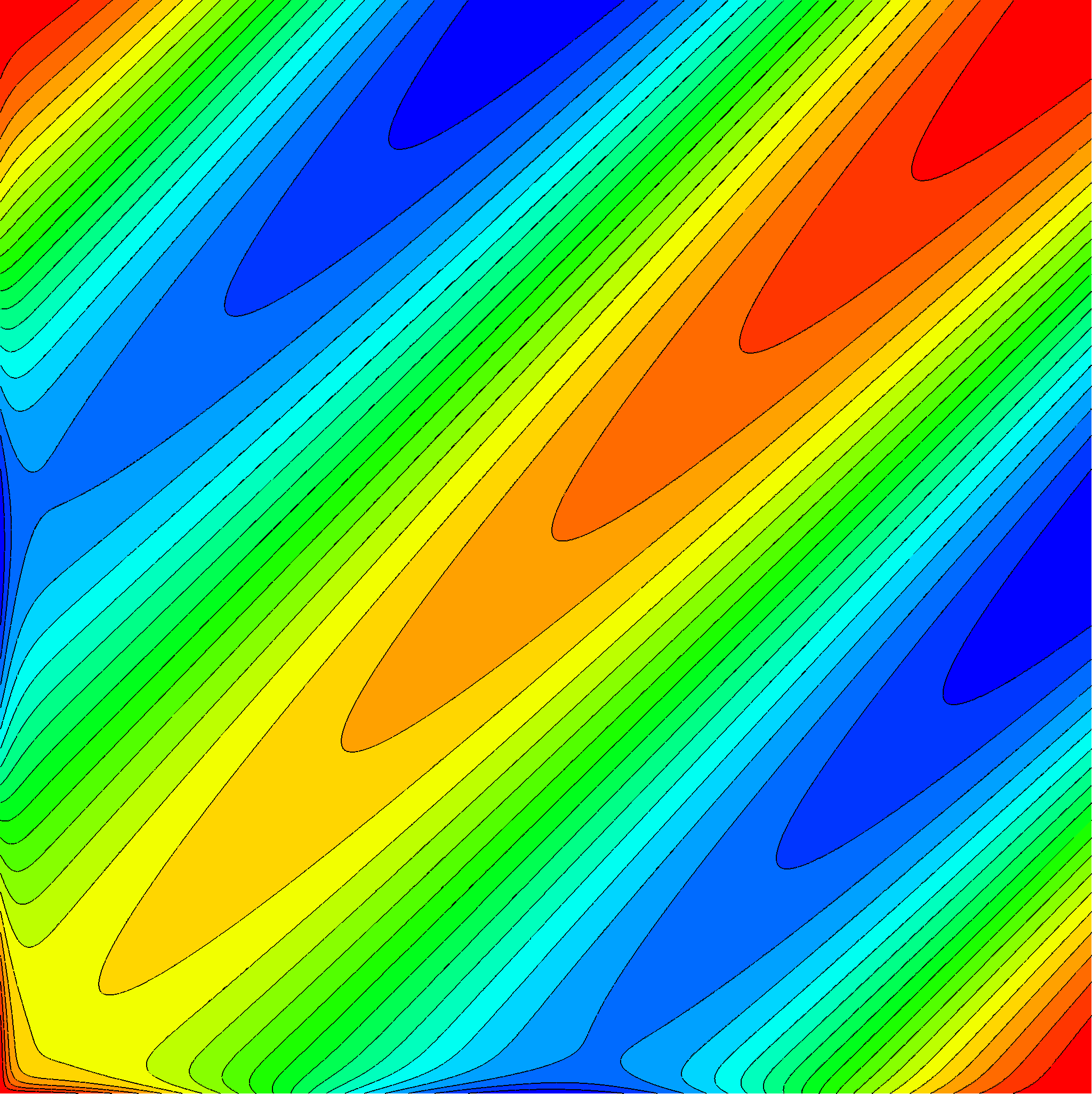}}
  \caption{Primal and dual solution for the convection-diffusion test case.}
\end{figure}

We consider a scalar convection-diffusion equation and apply a source function $s$  that allows us to prescribe the analytical solution in advance:
\begin{alignat*}{2}
\nabla\cdot(w,w)-\epsilon\Delta w&=s\qquad &&(x,y)\in\Omega=[0,1]^2,\\
w(x,y)&=0\qquad &&(x,y)\in\partial\Omega.
\end{alignat*}
We define $s:=s(x,y)$ such that
\begin{equation*}
w(x,y)=\left(x+\frac{e^{x/\epsilon}-1}{1-e^{1/\epsilon}}\right)\cdot\left(y+\frac{e^{y/\epsilon}-1}{1-e^{1/\epsilon}}\right)
\end{equation*}
is the solution to the equation (see Fig.~\ref{fig:bdrylayer_primal}, $\epsilon=0.01$). Depending on the value of $\epsilon$ a thin boundary layer evolves at the upper and right boundary.

The target functional of interest is the weighted total boundary flux, i.e.
\begin{equation*}
J(q,w)=\int_{\partial\Omega} \psi\left(w-\epsilon n\cdot q\right)\ds
\end{equation*}
where we choose $\psi(x,y)=\cos{\left(2\pi x\right)}\cos{\left(2\pi y\right)}$. Please note, that this weighting function represents the boundary conditions for the adjoint problem, i.e.
\begin{alignat*}{2}
-\nabla\cdot(\widetilde{w},\widetilde{w})-\epsilon\Delta \widetilde{w}&=0\qquad &&(x,y)\in\Omega=[0,1]^2\\
\widetilde{w}(x,y)&=\psi(x,y)\qquad &&(x,y)\in\partial\Omega.
\end{alignat*}
The adjoint solution can be seen in Fig.~\ref{fig:bdrylayer_dual}. One can clearly observe that the direction of advection is reversed compared to the primal solution. Furthermore, steep gradients can be observed in the lower left corner due to the Dirichlet boundary conditions.

The true value of the target functional of the manufactured solution depends only on $\epsilon$ and is given by 
\begin{equation*}
J_{\epsilon}(q,w)=-\frac{2\epsilon}{1+\left(2\pi\epsilon\right)^2}
\end{equation*}
so that $J_{0.01}(q,w)\approx -0.0199214$.

\begin{figure}[h]
\centering
  \subfloat[$p=1$]{\includegraphics{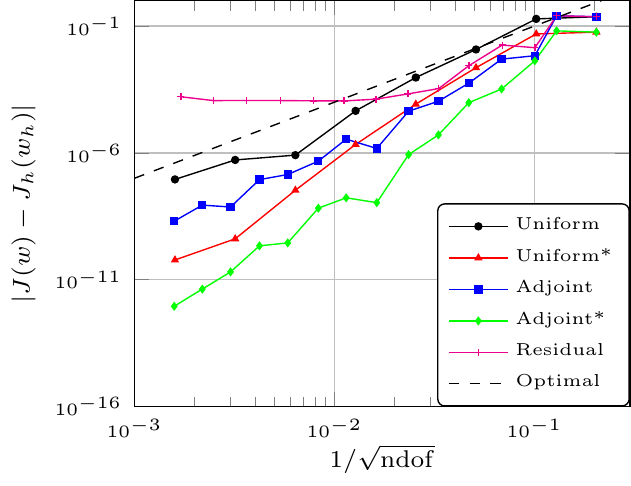}}\qquad
  \subfloat[$p=2$]{\includegraphics{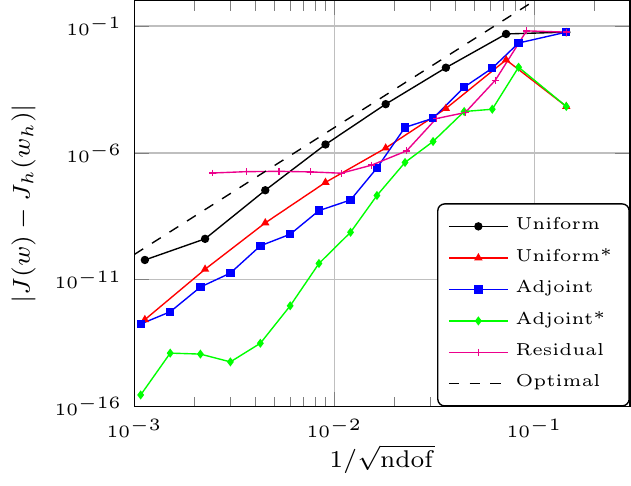}}\\
  \subfloat[$p=3$]{\includegraphics{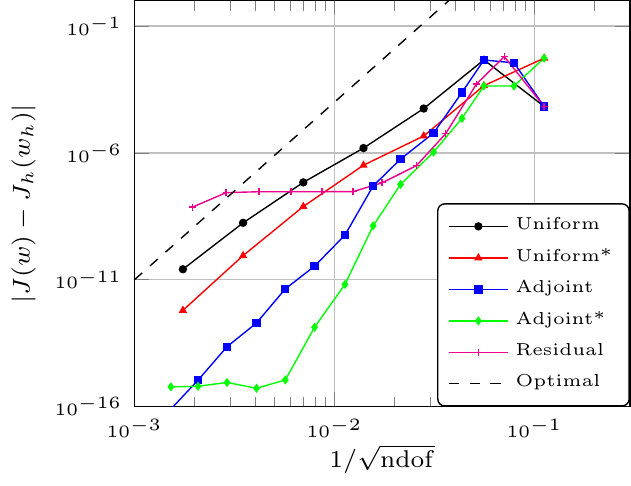}}\qquad
  \subfloat[$p=4$]{\includegraphics{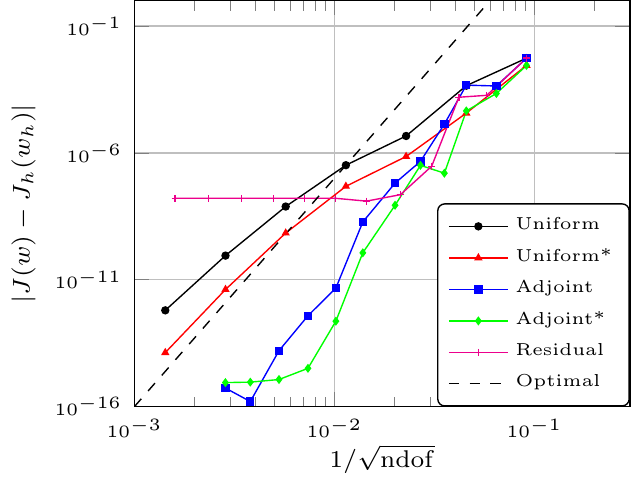}}
  \caption{Error with respect to degrees of freedom for the scalar boundary layer test case. A * denotes adjoint-based corrected values. Optimally, the error should converge as $\mathcal{O}\left(h^{2p+1}\right)$.}	
  \label{fig:bdrylayer_ndof}
\end{figure}

For uniform refinement the error rates (see Fig.\,\ref{fig:bdrylayer_ndof}) are approaching the optimal value of $2p+1$. By computing the adjoint even for uniform refinement we can correct the computed functional and obtain an error which is roughly two orders of magnitude smaller.

\begin{figure}[h]
\centering
  \subfloat[Residual-based refined mesh\label{fig:bdrylayer_mesh_h_res}]{
  \includegraphics[width=0.4\textwidth]{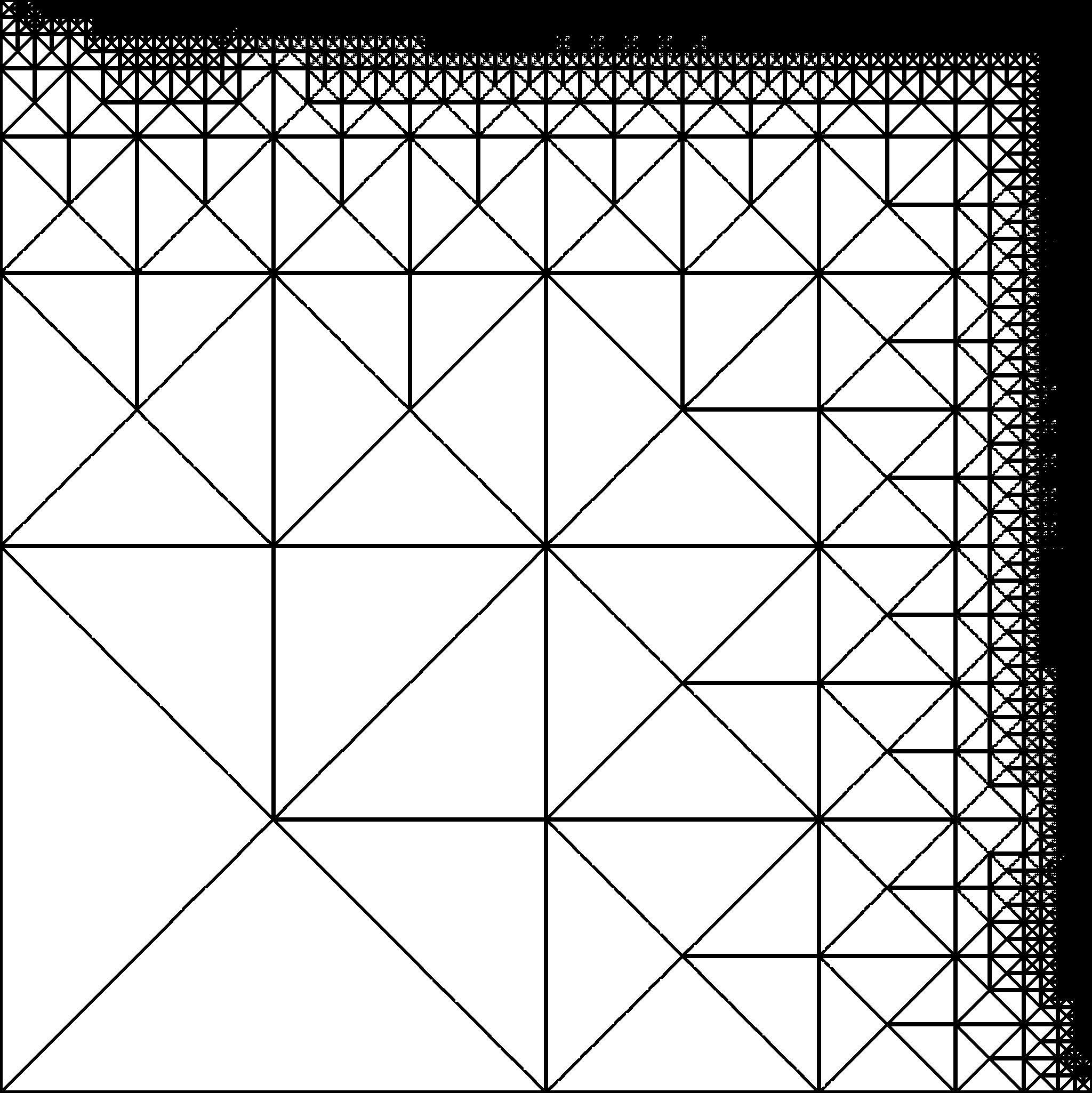}}\qquad
  \subfloat[Adjoint-based refined mesh\label{fig:bdrylayer_mesh_h}]{
  \includegraphics[width=0.4\textwidth]{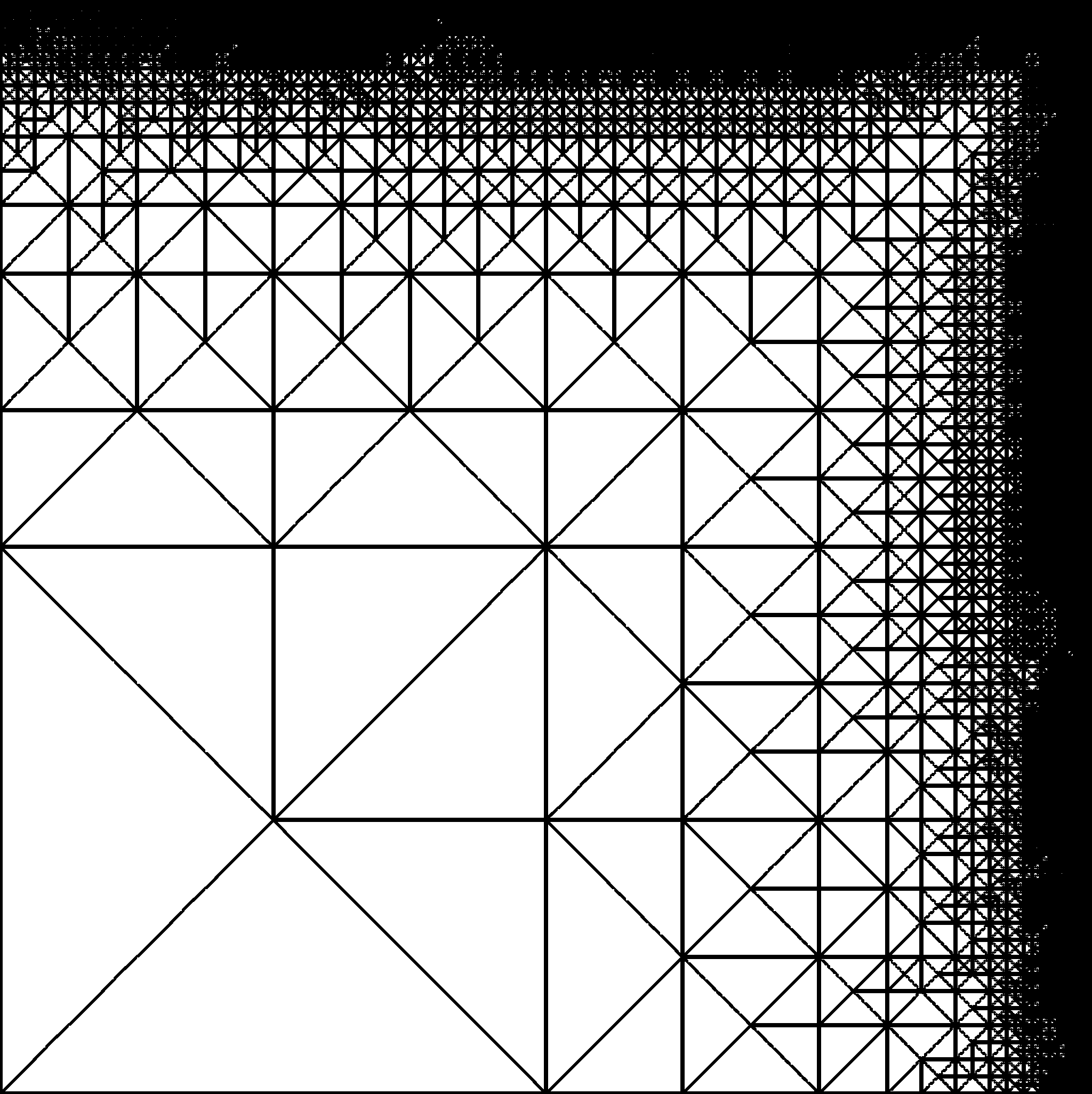}}
  \caption{Adapted meshes for the convection-diffusion test case ($p=2$, $n_e=18000$).}
\end{figure}

We use $\theta=0.15$ for both adjoint- and residual-based refinement. The latter shows a good behavior during the first adaptation steps, then however stalls. Adjoint-based refinement yields the lowest errors, especially when using the error estimate to correct the computed functional. In Fig.~\ref{fig:bdrylayer_mesh_h} the adjoint-based refined mesh can be seen. The sensor concentrates on the boundaries where the thin boundary layer is located. It, however, also refines regions which are upstream of the boundary layer which are important due to the strong convection. The residual-based indicator on the other hand refines only along the boundaries (see Fig.~\ref{fig:bdrylayer_mesh_h_res}).

\FloatBarrier
\subsection{Smooth Bump}
We consider subsonic, inviscid flow in a channel of height 0.8 and length 3 with an inlet Mach number of ${\rm Ma}_{\infty}=0.5$. The shape of the lower wall is a Gaussian of amplitude 0.0625 which has its center at mid-length of the channel. The flow is perfectly smooth, so that there should be no entropy production. This makes it a very good measure for mesh convergence. The target functional is given as follows
\begin{equation*}
J(w):=\|\Delta s\|_2:=\sqrt{\frac1{|\Omega|}\int_{\Omega}\left(\frac{p/\rho^{\gamma}-p_{\infty}/\rho^{\gamma}_{\infty}}{p_{\infty}/\rho^{\gamma}_{\infty}}\right)^2\dx}
\end{equation*}
which represents the non-dimensional entropy production.
 
\begin{figure}[h]
\centering
\subfloat[Baseline mesh\label{fig:bump_mesh_coarse}]{
\includegraphics[width=0.8\textwidth]{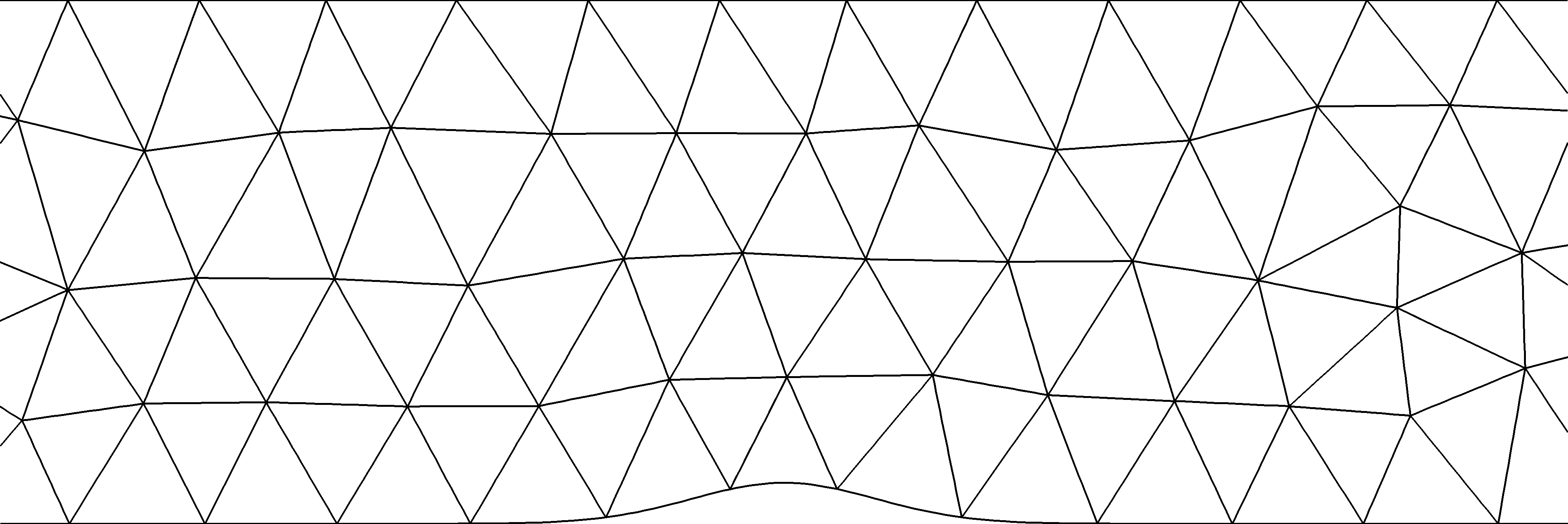}}\\
\subfloat[Residual-based refined mesh ($p=2$, $n_e\approx 18000$)\label{fig:bump_mesh_h_res}]{
\includegraphics[width=0.8\textwidth]{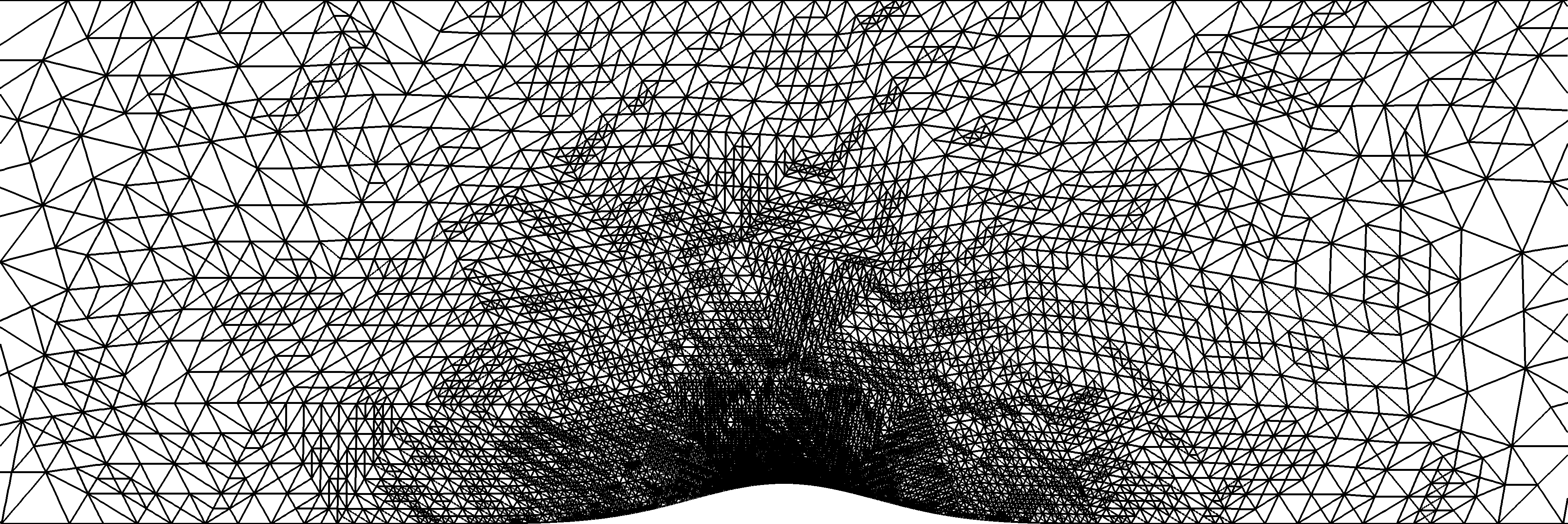}}\\
\subfloat[Adjoint-based refined mesh ($p=2$, $n_e\approx 18000$)\label{fig:bump_mesh_h}]{
\includegraphics[width=0.8\textwidth]{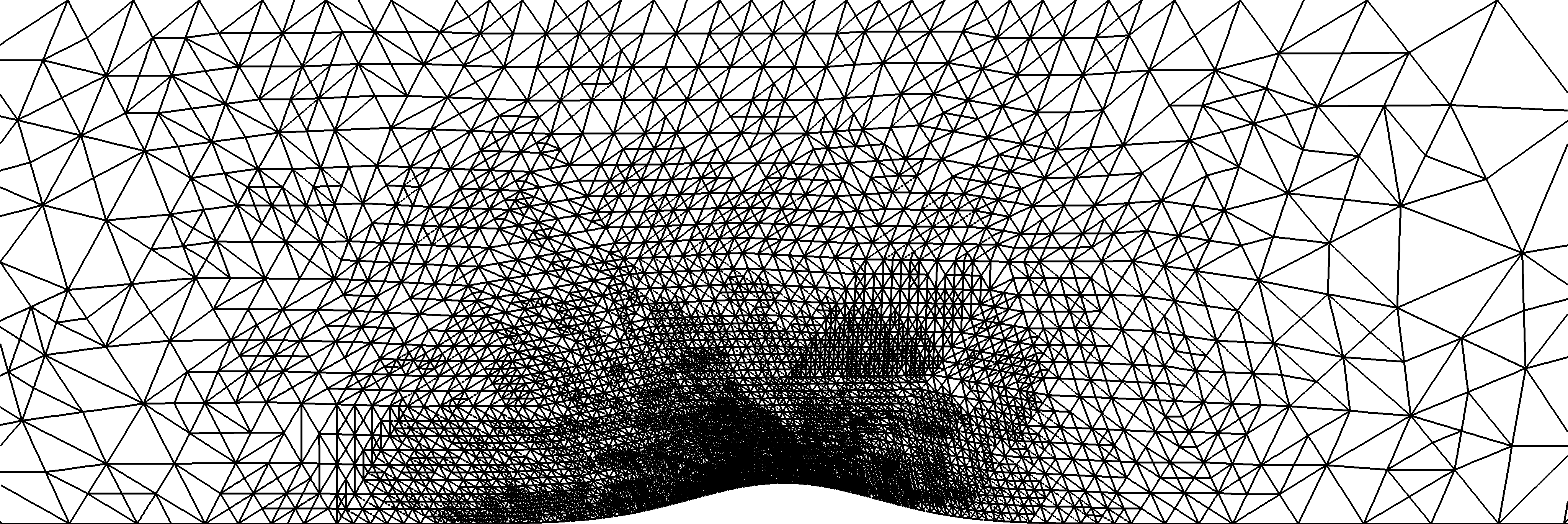}}
\caption{Meshes for the smooth bump test case (${\rm Ma}_{\infty}=0.5$, $\alpha=0^{\circ}$)}
\end{figure}

We apply Riemann invariant in- and outflow conditions at the left and right boundary and slip-wall conditions along the upper and lower boundary. The Mach number is set to $\rm Ma=0.5$. The baseline mesh, for both uniform and adaptive refinement, consists of 124 triangular elements (see Fig.~\ref{fig:bump_mesh_coarse}). The element transformation for those elements sharing an edge with this boundary are of 5th order. This is crucial in order to realize the theoretical order of convergence which is bounded by the minimum order of trial functions and element transformations at curved boundaries, respectively (see \citet{bassi1997dgeuler}).

\begin{figure}[h]
\centering
\subfloat[Mach number ($p=2$)\label{fig:bump_ma}]{
\includegraphics[width=0.8\textwidth]{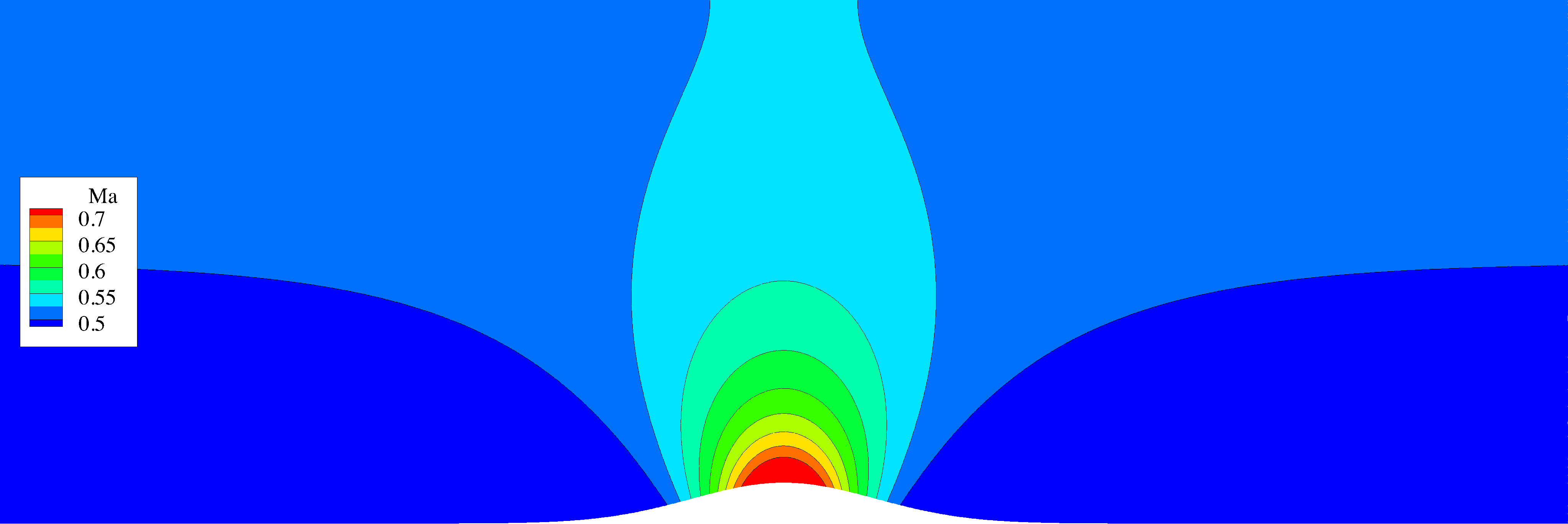}}\\
\subfloat[Adjoint x-momentum ($p=3$)\label{fig:bump_adjoint2}]{
\includegraphics[width=0.8\textwidth]{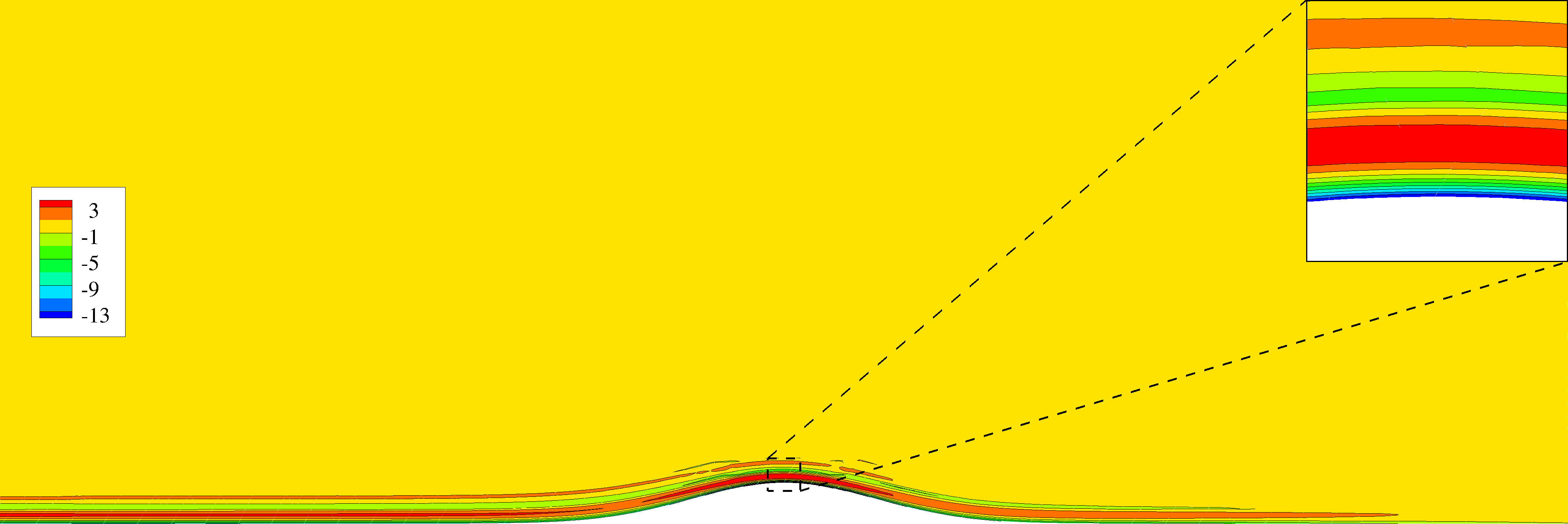}}
\caption{Primal and dual (entropy) solution for the smooth bump test case (${\rm Ma}_{\infty}=0.5$, $\alpha=0^{\circ}$).}
\end{figure}

In Fig.~\ref{fig:bump_ma} and \ref{fig:bump_adjoint2} contours of the Mach number and the adjoint x-momentum with respect to the entropy production are given. The primal solution is axisymmetric. The dual solution exhibits some interesting features. Please note, that the target functional is a volume integral and thus acts as a source term within the adjoint equation. The dual solution is close to zero in the majority of the domain. It attains higher (or lower) values towards the curved boundary where a thin layer develops which gets convected upstream (please recall that the flow direction is reversed in the adjoint equations). The extreme values are reached just at the top of the bump (see close-up in Fig.~\ref{fig:bump_adjoint2}).

\begin{figure}[h]
\centering
  \subfloat[$p=1$]{\includegraphics{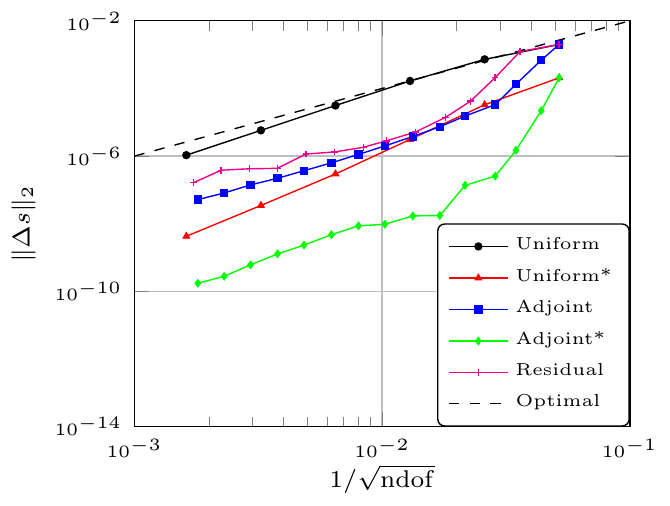}}\qquad
  \subfloat[$p=2$]{\includegraphics{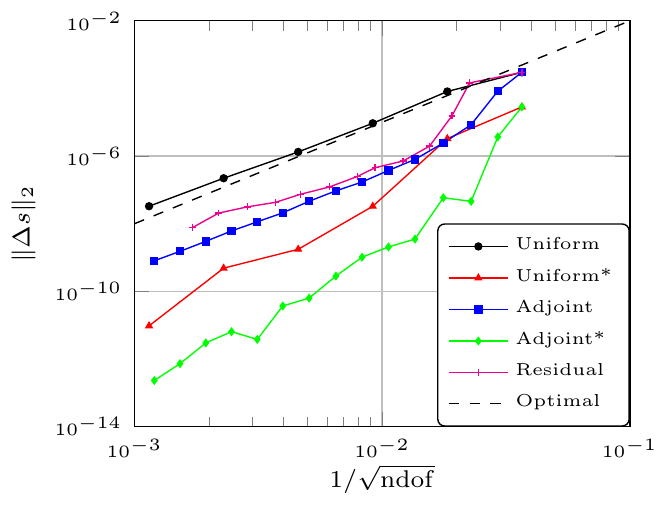}}\\
  \subfloat[$p=3$]{\includegraphics{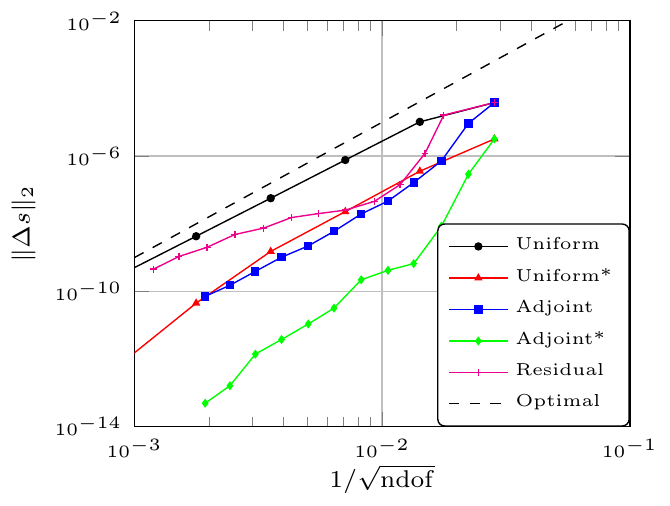}}\qquad
  \subfloat[$p=4$]{\includegraphics{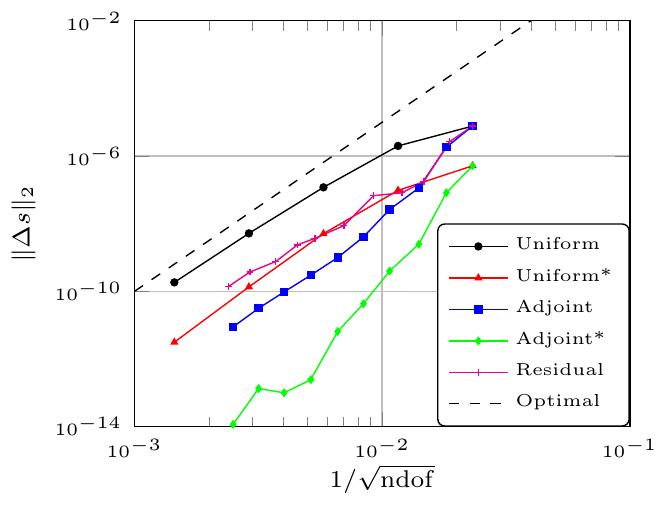}}
	\caption{Entropy error with respect to degrees of freedom. A * denotes adjoint-based corrected values. Optimally, the error should converge as $\mathcal{O}\left(h^{p+1}\right)$.}
	\label{fig:bump_ndof}
\end{figure}

Under uniform refinement we obtain the expected order of convergence, namely $\|\Delta s\|_2=\mathcal{O}\left(h^{p+1}\right)$ (see Fig.\,\ref{fig:bump_ndof}). Furthermore, it is important to note that a higher polynomial order for the ansatz functions is more efficient than a lower order. Additionally, we computed the adjoint solution on the uniformly refined meshes and corrected the computed error with the given estimate. One can see that the corrected error is approximately two orders of magnitude smaller than the actually computed error.

We choose $\theta=0.1$ for the adaptive runs. Both residual- and adjoint-based adaptation show better error convergence compared to uniform refinement. In Fig.~\ref{fig:bump_mesh_h_res} and Fig.~\ref{fig:bump_mesh_h} one can see that the region around the bump exhibits the major part of refinement. In the latter, refinement is however more confined to the bump. After a drastic initial drop, the adjoint-based refinement attains the optimal order of convergence whereas the residual-based refinement shows a more chaotic behavior. Here, the superiority of the adjoint-based error sensor becomes apparent. This advantage is even more amplified when considering the corrected error which is two to four orders of magnitude lower compared to the error attained by uniform refinement.

\FloatBarrier
\subsection{NACA~0012}

In the following, we are concerned with three test cases of compressible flow around the NACA~0012 airfoil: an inviscid subsonic, a transonic and a viscous subsonic. The NACA~0012 airfoil is defined by
\begin{equation*}
y=\pm 0.6\left(0.2969\sqrt{x}-0.1260x-0.3516x^2+0.2843x^3-0.1015x^4\right)
\end{equation*}
with $x\in [0,1]$. Using this definition, the airfoil would have a finite trailing edge thickness of .252~\%. In order to obtain a sharp trailing edge we modify the $x^4$ coefficient, i.e.
\begin{equation*}
y=\pm 0.6\left(0.2969\sqrt{x}-0.1260x-0.3516x^2+0.2843x^3-0.1036x^4\right).
\end{equation*}

We use 4th-order curved elements along the boundaries and employ a circular far field boundary which is located more than 1000 chord lengths away from the airfoil.

The functional of interest can be the drag or the lift coefficient along the wall-boundary $\partial\Omega_w$, and is given by
\begin{equation*}
J(w)=\int_{\partial\Omega_w}p(w)\beta\cdot n\ds
\end{equation*}
with $\beta_d=\frac1{C_{\infty}}(\cos\alpha,\sin\alpha)^T$ for drag and $\beta_l=\frac1{C_{\infty}}(-\sin\alpha,\cos\alpha)^T$ for lift. $\alpha$ denotes the angle of attack and $C_{\infty}$ is a normalized reference value given by 
$C_{\infty}=\frac12\gamma p_{\infty}\mathrm Ma^2_{\infty}l$ where $l$ is the chord length of the airfoil, $p_{\infty}$ and ${\rm Ma}_{\infty}$ are the free-stream pressure and Mach number. If viscous effects are taken into account the drag and lift coefficients include the shear friction forces and are given by
\begin{equation*}
J(q,w)=\int_{\partial\Omega_w}p(w)\beta\cdot n-(\tau(w,q)\beta)\cdot n\ds.
\end{equation*}

In order to render the scheme adjoint consistent we have to evaluate the target functional with the boundary state $w_{\partial\Omega}(w)$ and the boundary flux $f_{v, \partial \Omega}$.

\FloatBarrier
\subsubsection{Inviscid, Subsonic Flow}

We consider inviscid, subsonic flow around the NACA0012 airfoil with a free-stream Mach number of $\rm Ma_{\infty}=0.5$ and an angle of attack $\alpha=2^{\circ}$.

This flow has one very distinct feature, namely the singularity at the sharp trailing edge. Furthermore, the solution exhibits high gradients in the leading edge region due to the stagnation of the flow.

The adaptive routine is driven by an error estimate for the drag coefficient. In each step, the topmost 2\% elements with respect to error contribution are marked for refinement ($\theta=0.02$). For higher $\theta$ too many elements contributing only little to the overall error were marked during the first refinement steps. In order to compute the error in the drag coefficient, a reference value was obtained on an adjoint-based adapted mesh with approximately $2.6\cdot 10^5$ degrees of freedom.

\begin{figure}[h]
\centering
\subfloat[Mach number ($p=2$)\label{fig:naca_ma05_al2_ma}]{
\includegraphics[width=0.8\textwidth]{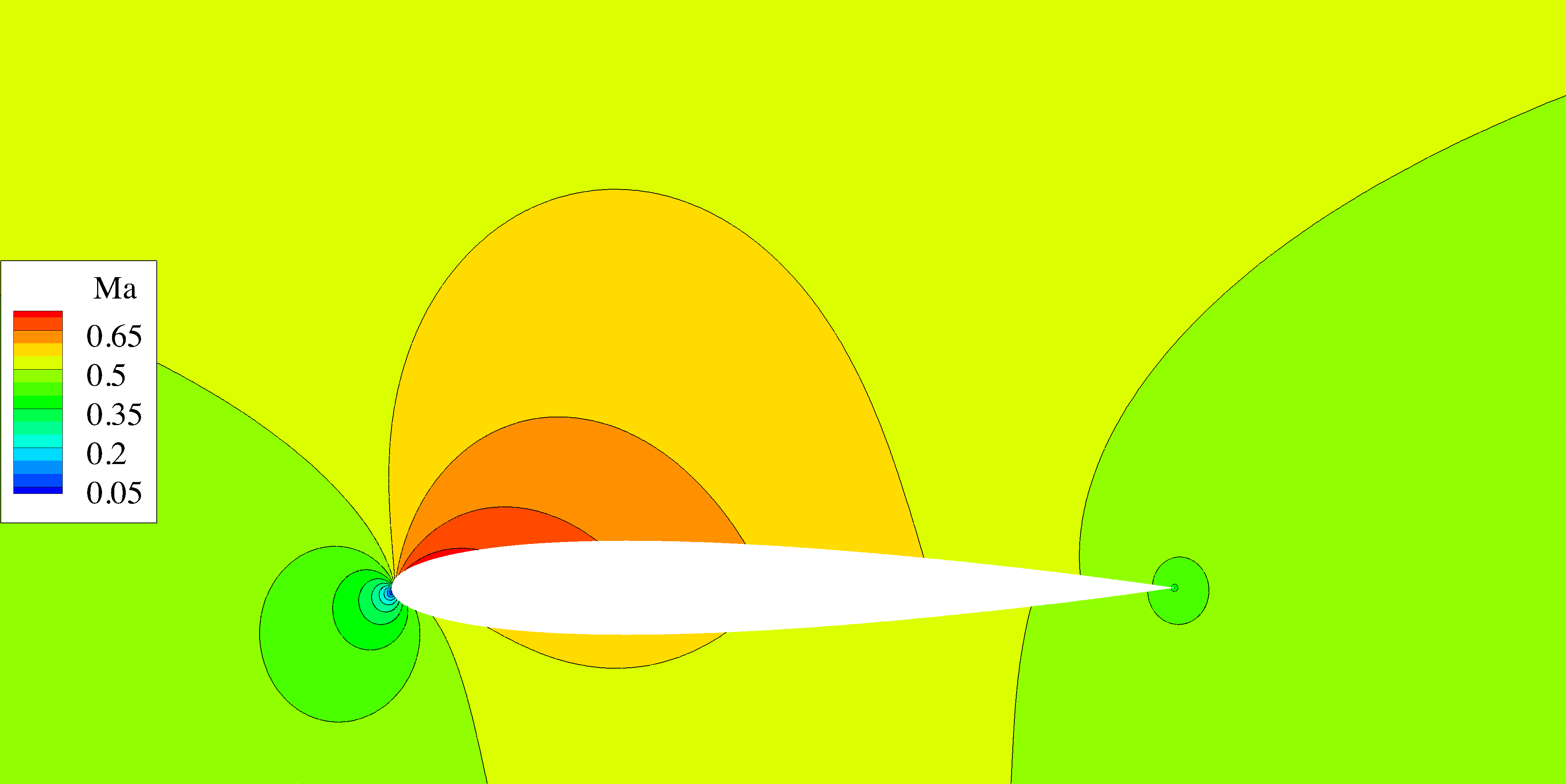}}\\
\subfloat[Adjoint x-momentum ($p=3$)\label{fig:naca_ma05_al2_adjoint2}]{
\includegraphics[width=0.8\textwidth]{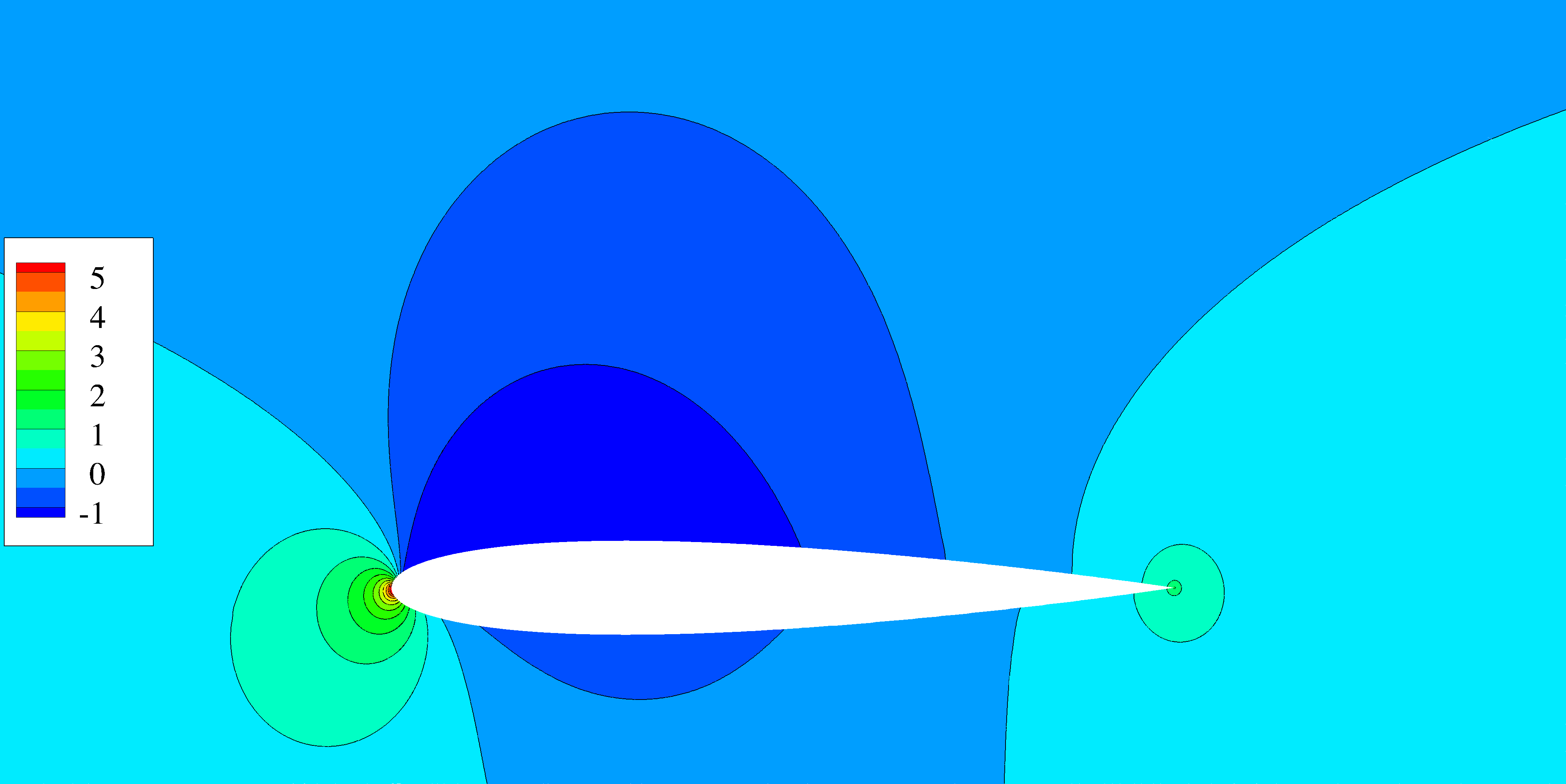}}
\caption{Primal and dual (drag) solution for the subsonic Euler test case (${\rm Ma}_{\infty}=0.5$, $\alpha=2^{\circ}$).}
\end{figure}

In Fig.~\ref{fig:naca_ma05_al2_ma} and \ref{fig:naca_ma05_al2_adjoint2} contours of the Mach number and the adjoint x-momentum with respect to drag are given. The approximate adjoint solution is very smooth due to the adjoint-consistent treatment of boundary conditions and target functionals. One can see that the adjoint solution shows peaks at leading and trailing edge. This indicates the importance of these two locations for the drag coefficient.

\begin{figure}[h]
\centering
  \subfloat[$p=1$]{\includegraphics{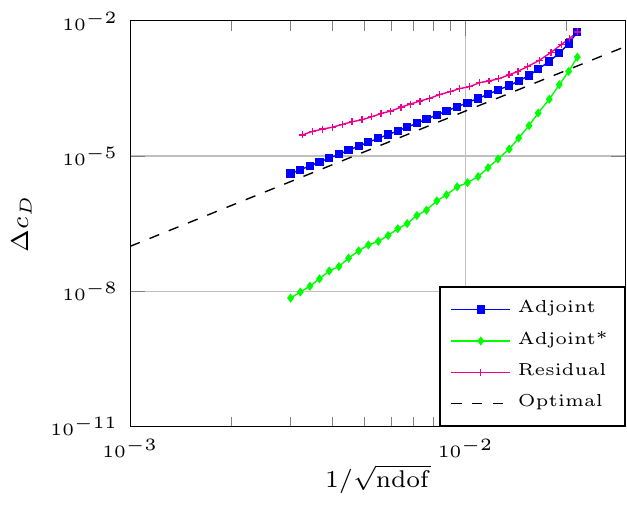}}\qquad
  \subfloat[$p=2$]{\includegraphics{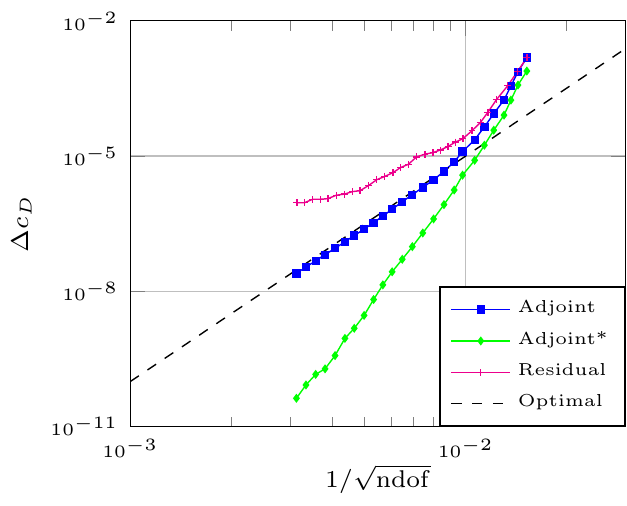}}\\
  \subfloat[$p=3$]{\includegraphics{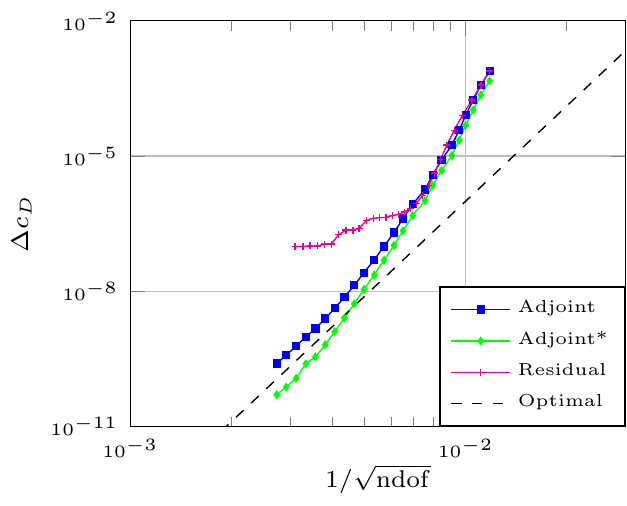}}\qquad
  \subfloat[$p=4$]{\includegraphics{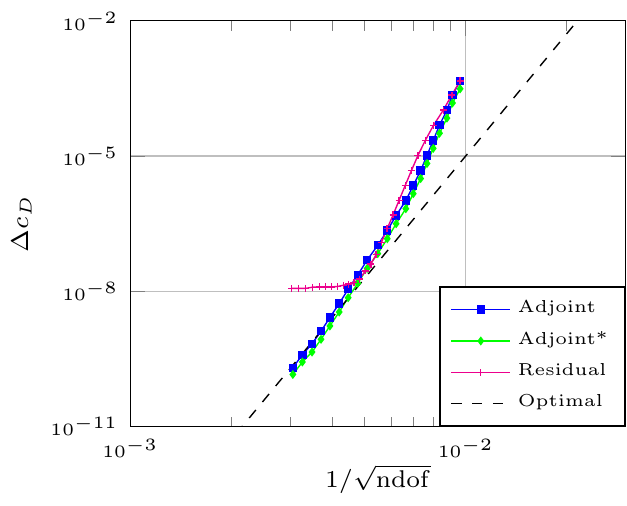}}  
  \caption{Drag convergence with respect to degrees of freedom (${\rm Ma}_{\infty}=0.5$, $\alpha=2^{\circ}$). A * denotes adjoint-based corrected values. Optimally, the error should converge as $\mathcal{O}\left(h^{2p+1}\right)$.}
  \label{fig:naca_ma05_al2_ndof}
\end{figure}

\begin{figure}[h]
\centering
\subfloat[Baseline mesh\label{fig:naca_ma05_al2_mesh_coarse}]{
\includegraphics[width=0.8\textwidth]{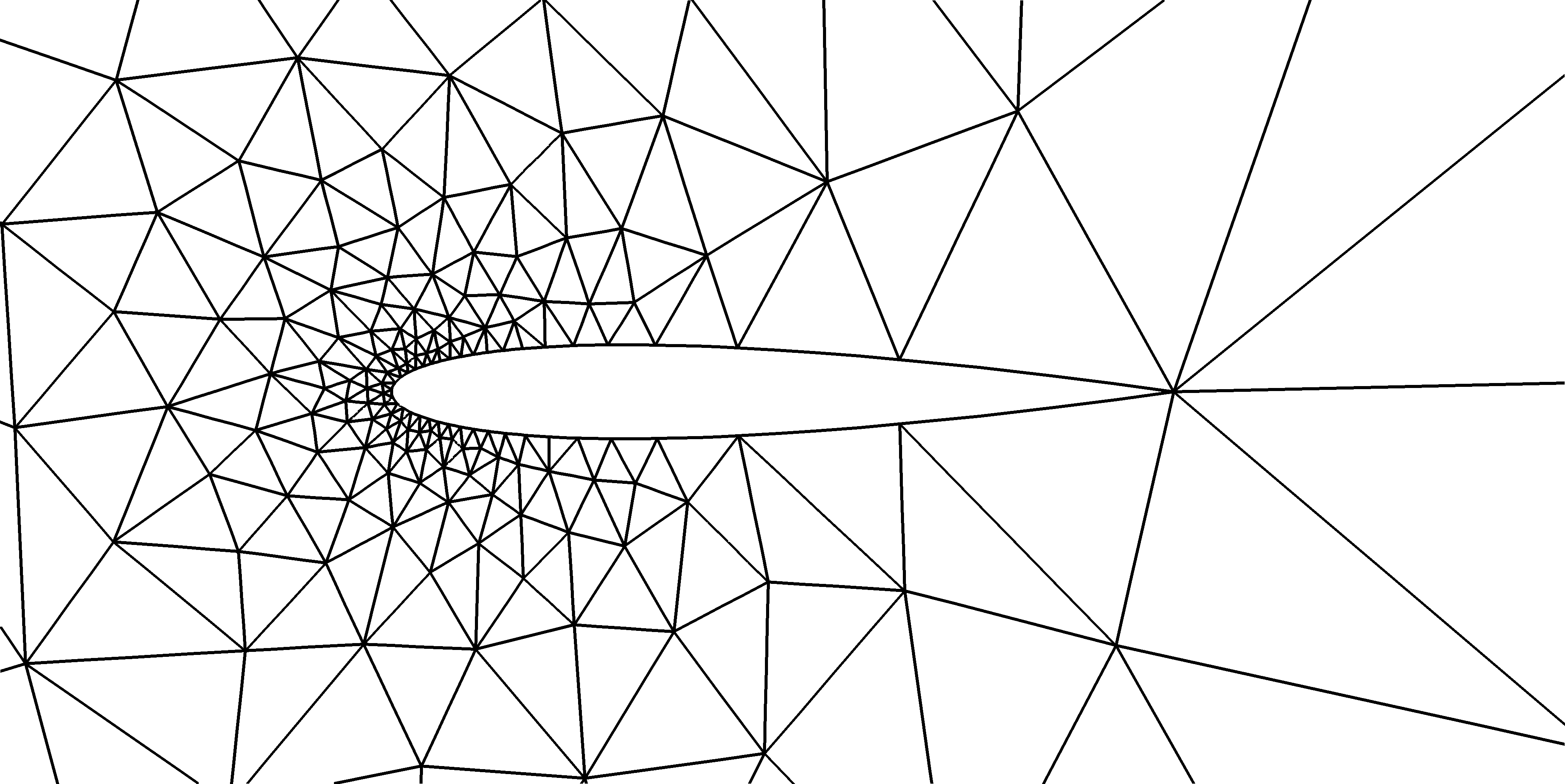}}\\
\subfloat[Residual-based refined mesh ($p=2$, $n_e\approx 18000$)\label{fig:naca_ma05_al2_mesh_h_res}]{
\includegraphics[width=0.8\textwidth]{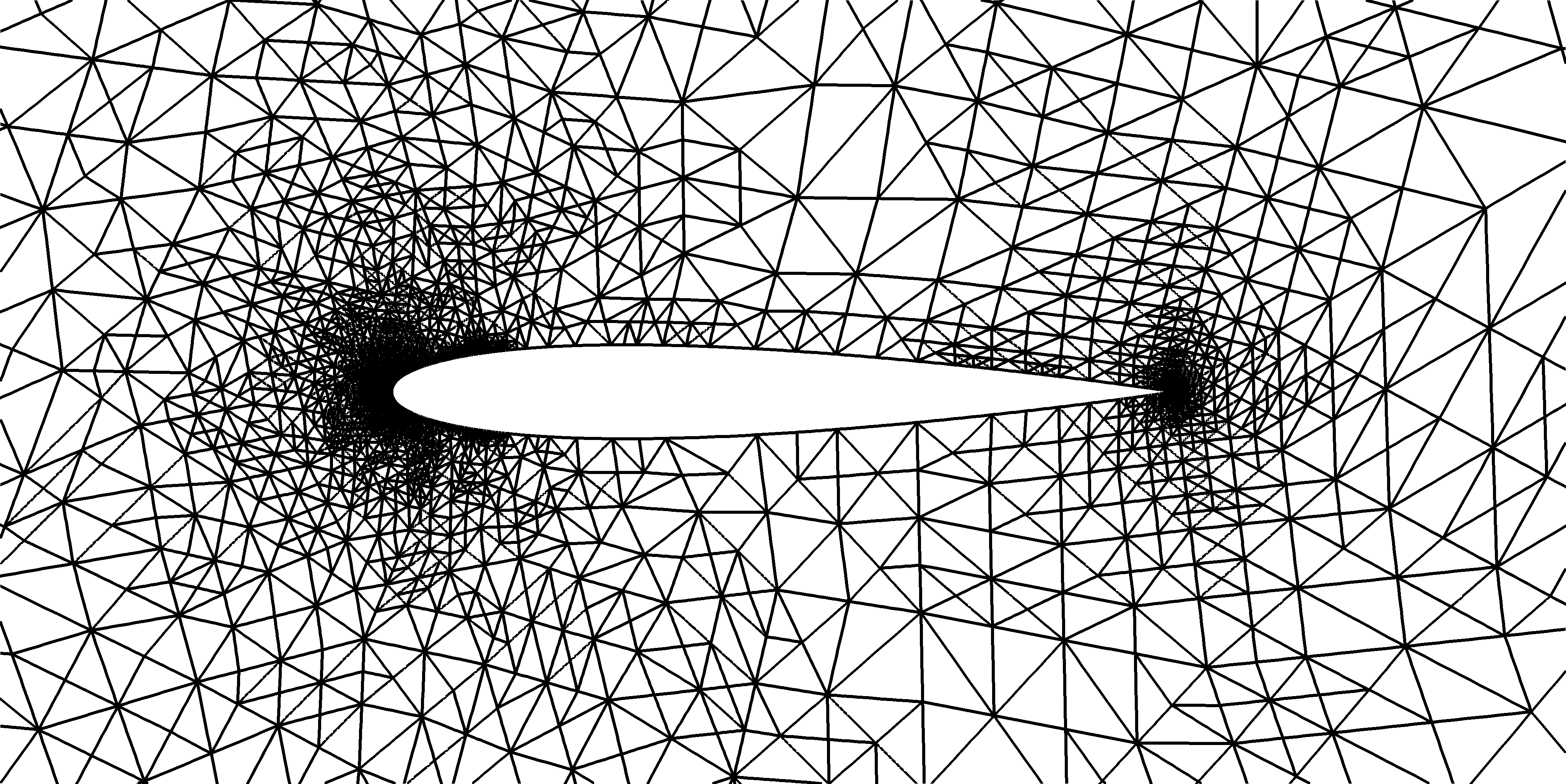}}\\
\subfloat[Adjoint-based refined mesh ($p=2$, $n_e\approx 18000$)\label{fig:naca_ma05_al2_mesh_h}]{
\includegraphics[width=0.8\textwidth]{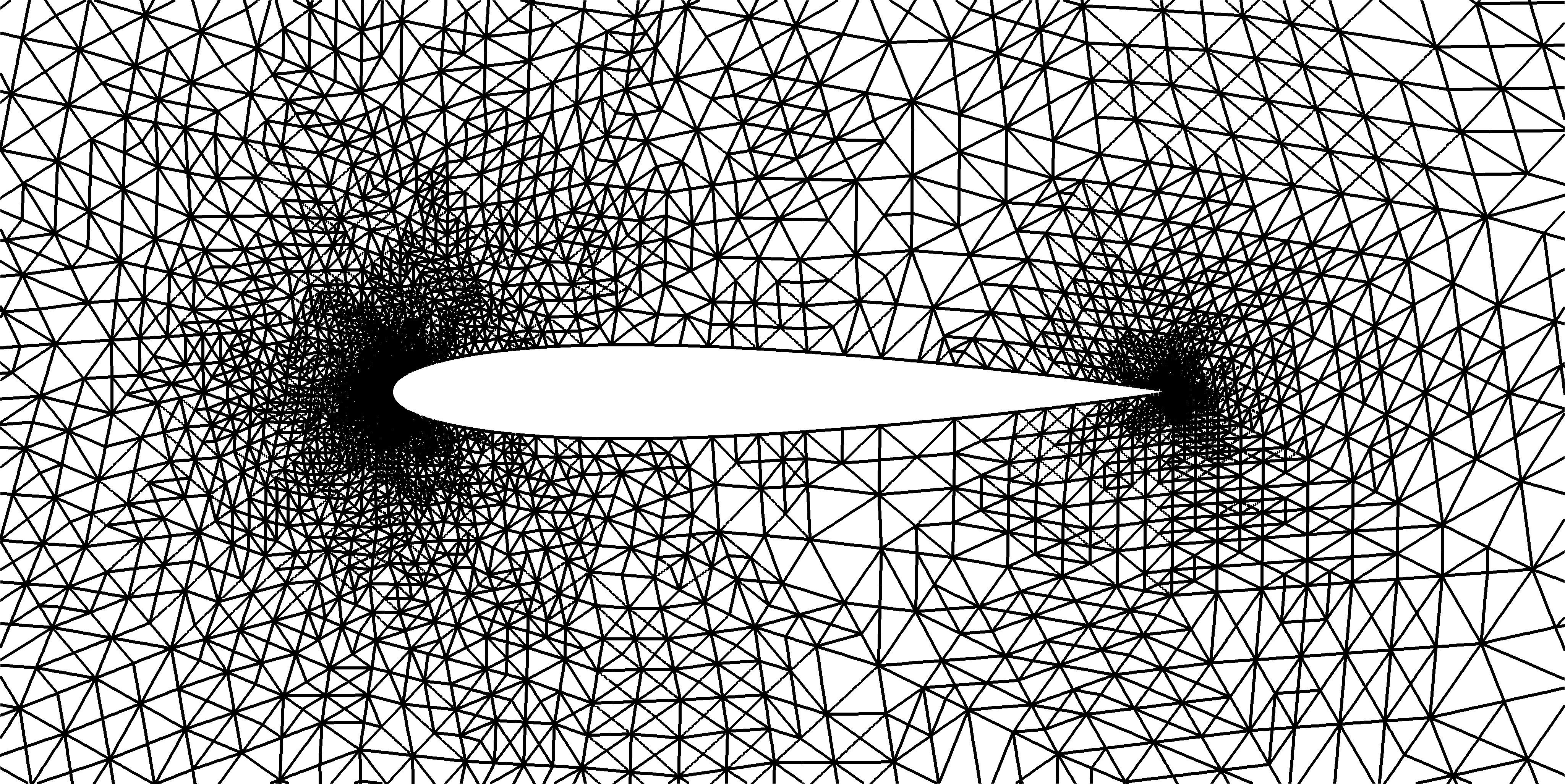}}
\caption{Meshes for the subsonic Euler test case (${\rm Ma}_{\infty}=0.5$, $\alpha=2^{\circ}$)}
\end{figure}

\begin{figure}
\centering
\includegraphics{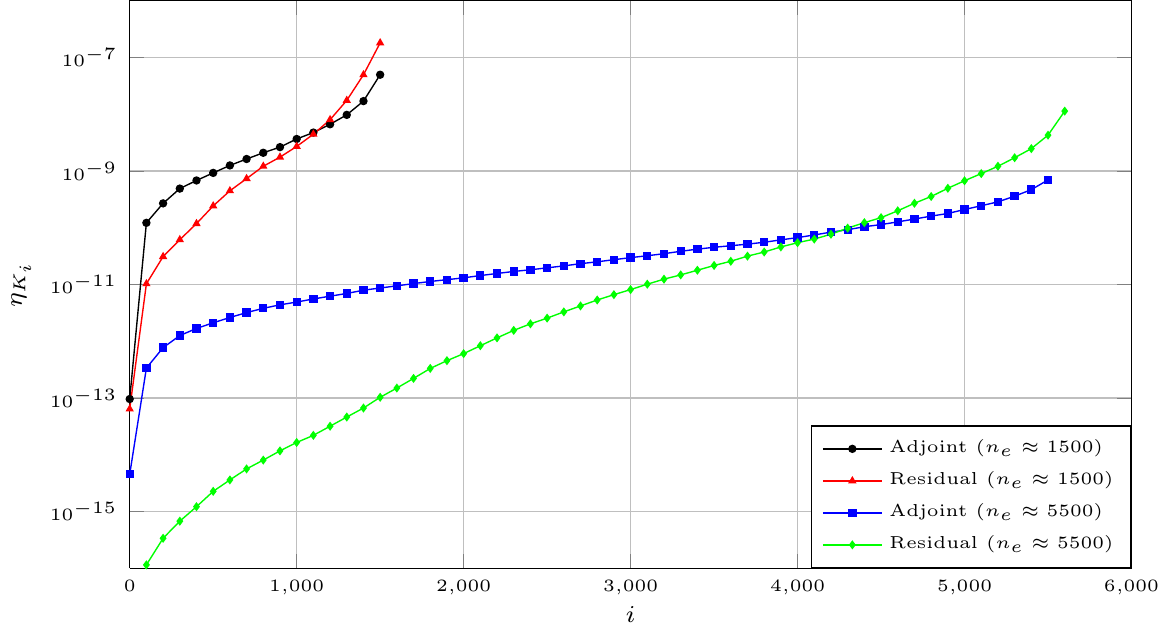}
\caption{Element-wise adjoint-based error estimate (sorted by magnitude, see Sec.~\ref{sec:marking}) for adjoint- and residual-based refinement (every 100th point is plotted)}
\label{fig:naca_ma05_al2_errest}
\end{figure}

During the first refinement steps, both adjoint-based and residual-based adaptation perform very well (see Fig.~\ref{fig:naca_ma05_al2_ndof}). Both sense and refine the above mentioned regions which exhibit a large gradient and a singularity, respectively. As soon as these regions have a resolution which yields a local contribution to the error in drag comparable to other regions, the convergence path differs. The adjoint-based adaptation achieves the theoretical order of $2p+1$ whereas the error for the residual-based adaptation flattens and finally stagnates. By taking a look at the refined meshes the reason for these different behaviors becomes clear (see. Fig.~\ref{fig:naca_ma05_al2_mesh_h_res} and \ref{fig:naca_ma05_al2_mesh_h}). Even when both leading and trailing edge have enough resolution, i.e.\, the local contributions to the error in drag are homogeneously distributed throughout the mesh, the residual-based adaptation procedure keeps on refining these two regions. The adjoint-based indicator, on the other hand, begins to mark other elements in the neighborhood of the airfoil as well. In Fig.~\ref{fig:naca_ma05_al2_errest} we have plotted the adjoint-based error estimate at two refinement stages ($n_e\approx 1500$ and $n_e\approx 5500$) for both local adjoint- and residual-based refinement (please note, that we computed the adjoint-based error estimate for the residual-based refined meshes afterwards). The element-wise error for the residual-based refined mesh spans a considerably larger range compared to the adjoint-based adapted mesh. Furthermore, the maximum error is bigger for residual-based refinement.

Using the adjoint-based error estimate to correct the drag value seems to be very efficient for $p=1,2$; for $p=3,4$ however, the difference between corrected and uncorrected values becomes very small. 

\pgfplotsset{cycle list={{blue, mark=square*,mark size=1}, {green,mark=diamond*,mark size=1},{magenta, mark=+,mark size=1}, {black,dashed}, {black}}}

\FloatBarrier
\subsubsection{Inviscid, Transonic Flow}

We consider inviscid, transonic flow around the NACA0012 airfoil with a free-stream Mach number of $\rm Ma_{\infty}=0.8$ and an angle of attack $\alpha=1.25^{\circ}$. This flow is dominated by a strong shock on the upper side of the airfoil and a weak shock on the lower side. Furthermore, the singularity at the trailing edge is present due to the slip-condition applied along the walls. We choose $\theta=0.1$ for the adaptive runs. In order to compute the error in the drag coefficient, a reference value  was obtained on an adjoint-based adapted mesh with approximately $2.3\cdot 10^5$ degrees of freedom.

\begin{figure}[h]
\centering
\subfloat[Mach number ($p=2$)\label{fig:naca_ma08_al125_ma}]{
\includegraphics[width=0.8\textwidth]{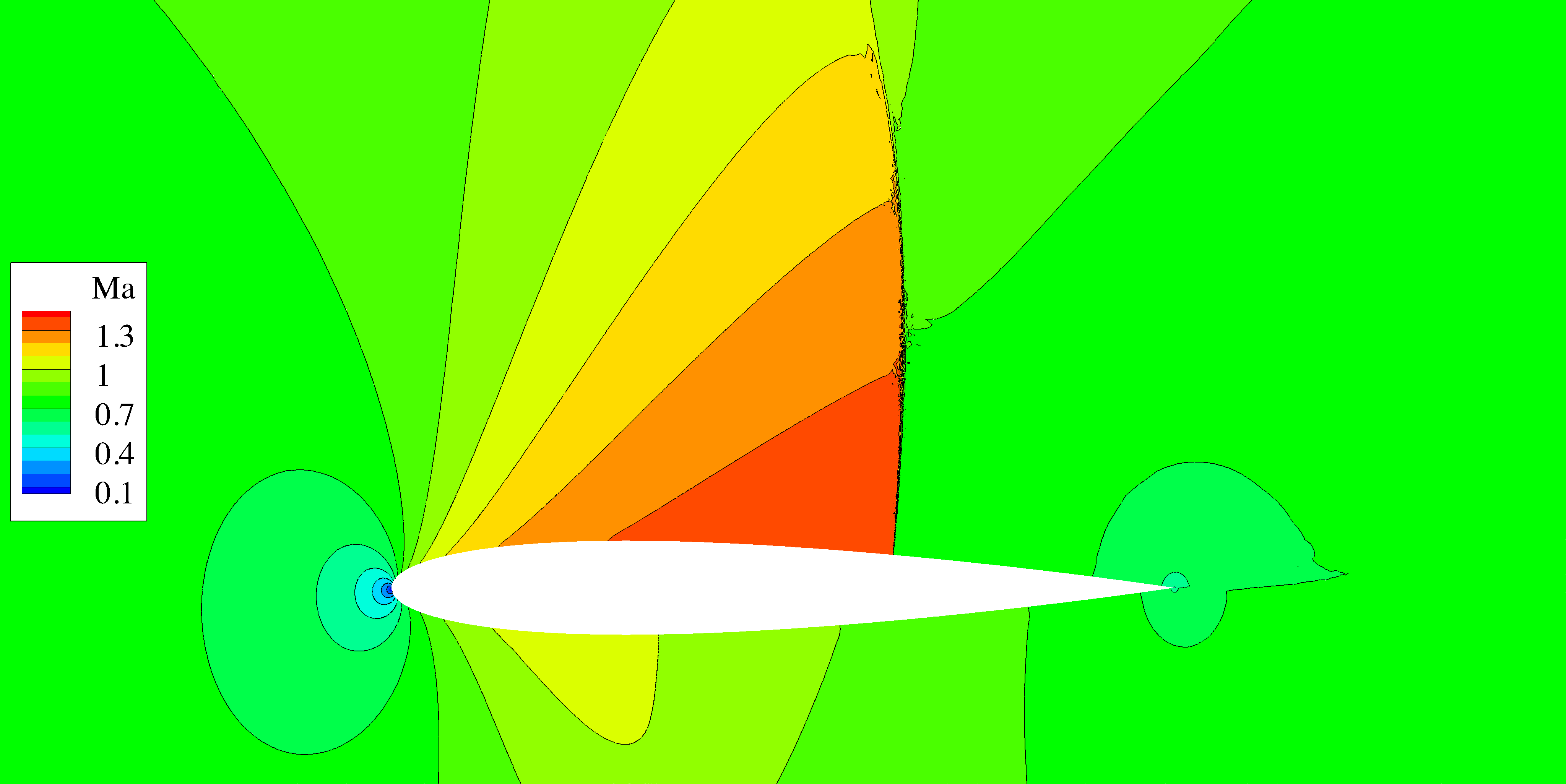}}\\
\subfloat[Adjoint x-momentum ($p=3$)\label{fig:naca_ma08_al125_adjoint2}]{
\includegraphics[width=0.8\textwidth]{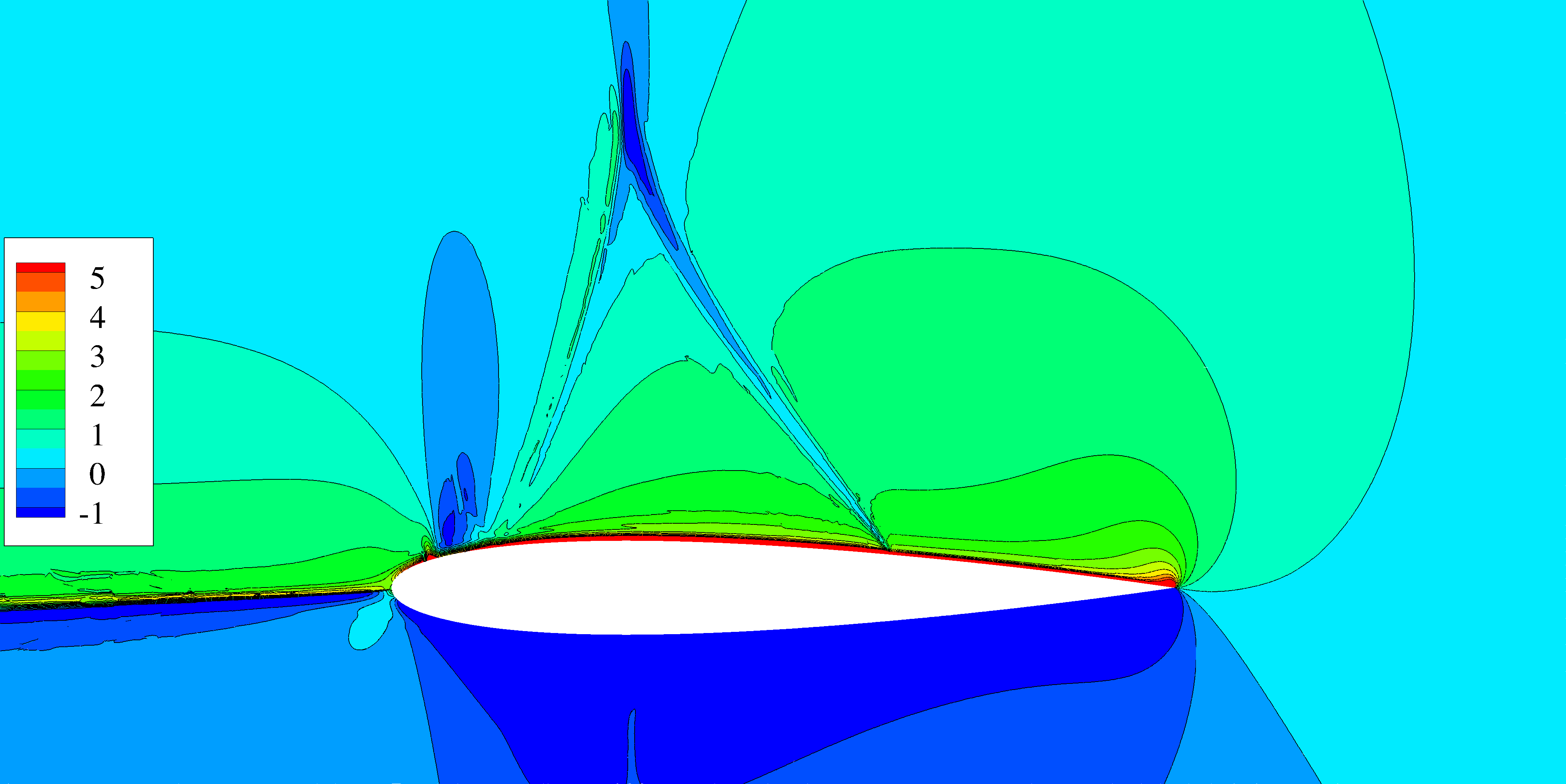}}
\caption{Primal and dual (drag) solution for the transonic Euler test case (${\rm Ma}_{\infty}=0.8$, $\alpha=1.25^{\circ}$).}
\end{figure}

In Fig.~\ref{fig:naca_ma08_al125_ma} and \ref{fig:naca_ma08_al125_adjoint2} contours of the Mach number and the adjoint x-momentum with respect to drag are given. Neither the primal nor the dual solution are smooth; both exhibit discontinuities. These, however, occur in different locations. Please note that no artificial viscosity is added in non-smooth regions of the adjoint solution.

\begin{figure}[h]
\centering
  \subfloat[$p=1$]{\includegraphics{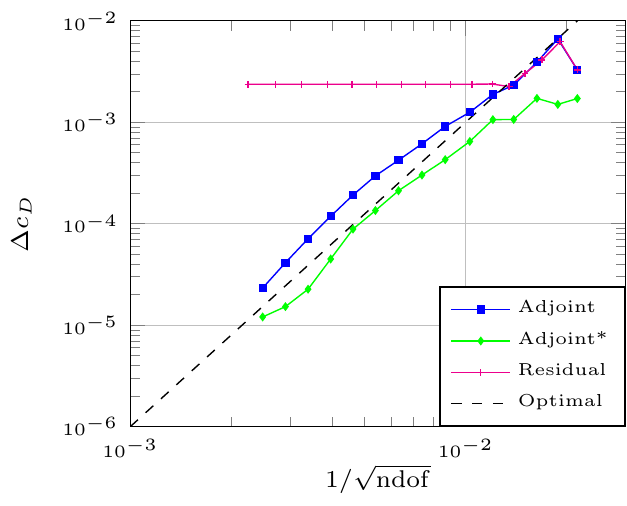}}\qquad
  \subfloat[$p=2$]{\includegraphics{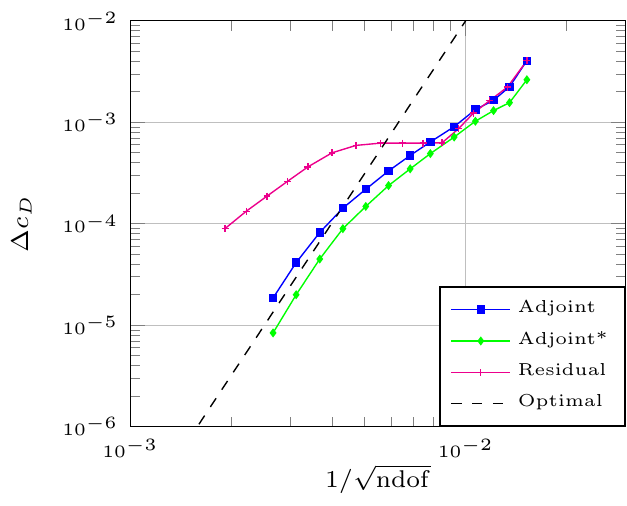}}\\
  \subfloat[$p=3$]{\includegraphics{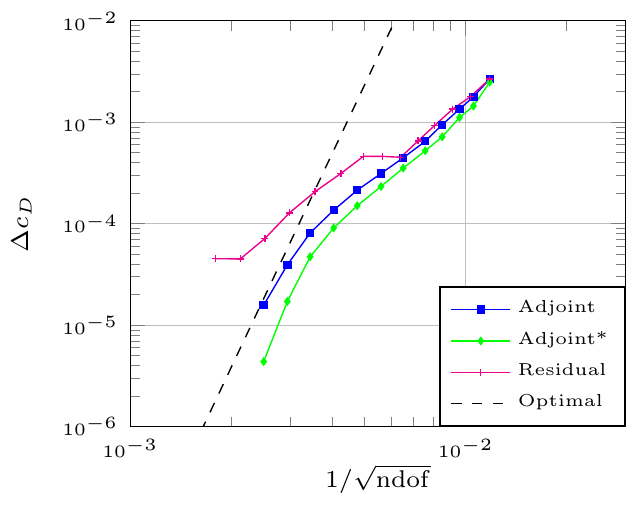}}\qquad
  \subfloat[$p=4$]{\includegraphics{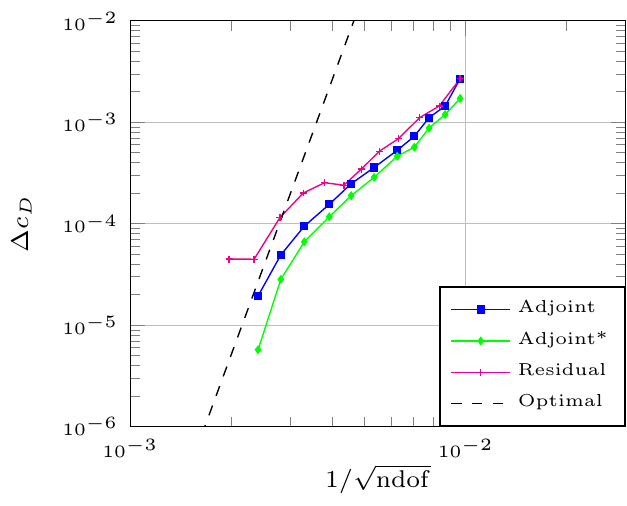}}
  \caption{Drag convergence with respect to degrees of freedom (${\rm Ma}_{\infty}=0.8$, $\alpha=1.25^{\circ}$). A * denotes adjoint-based corrected values. Optimally, the error should converge as $\mathcal{O}\left(h^{2p+1}\right)$.}
  \label{fig:naca_ma08_al125_ndof}
\end{figure}

\begin{figure}[h]
\centering
\subfloat[Baseline mesh\label{fig:naca_ma08_al125_mesh_coarse}]{
\includegraphics[width=0.8\textwidth]{naca_mesh_coarse.png}}\\
\subfloat[Residual-based refined mesh ($p=2$, $n_e\approx 18000$)\label{fig:naca_ma08_al125_mesh_h_res}]{
\includegraphics[width=0.8\textwidth]{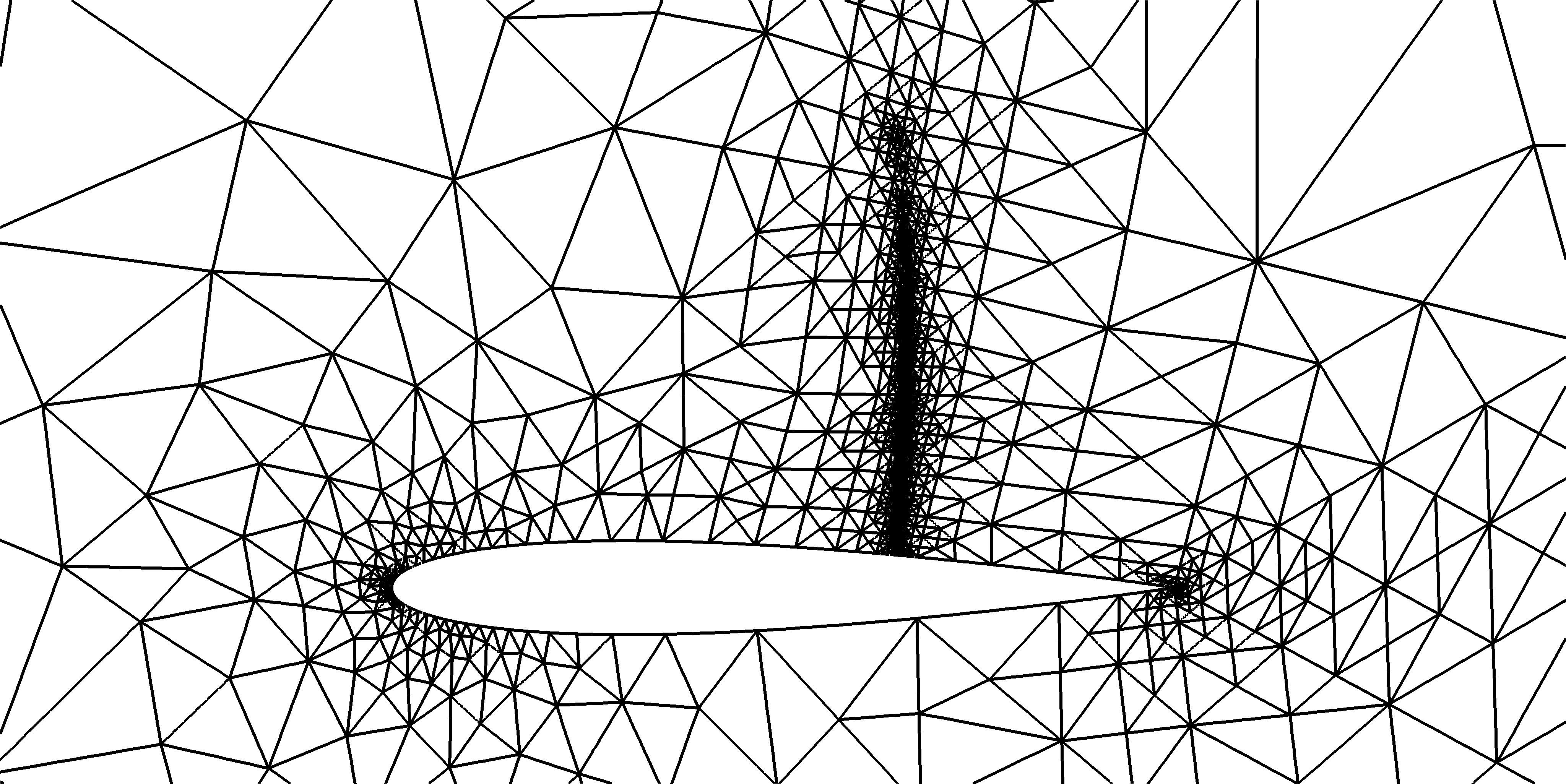}}\\
\subfloat[Adjoint-based refined mesh ($p=2$, $n_e\approx 18000$)\label{fig:naca_ma08_al125_mesh_h}]{
\includegraphics[width=0.8\textwidth]{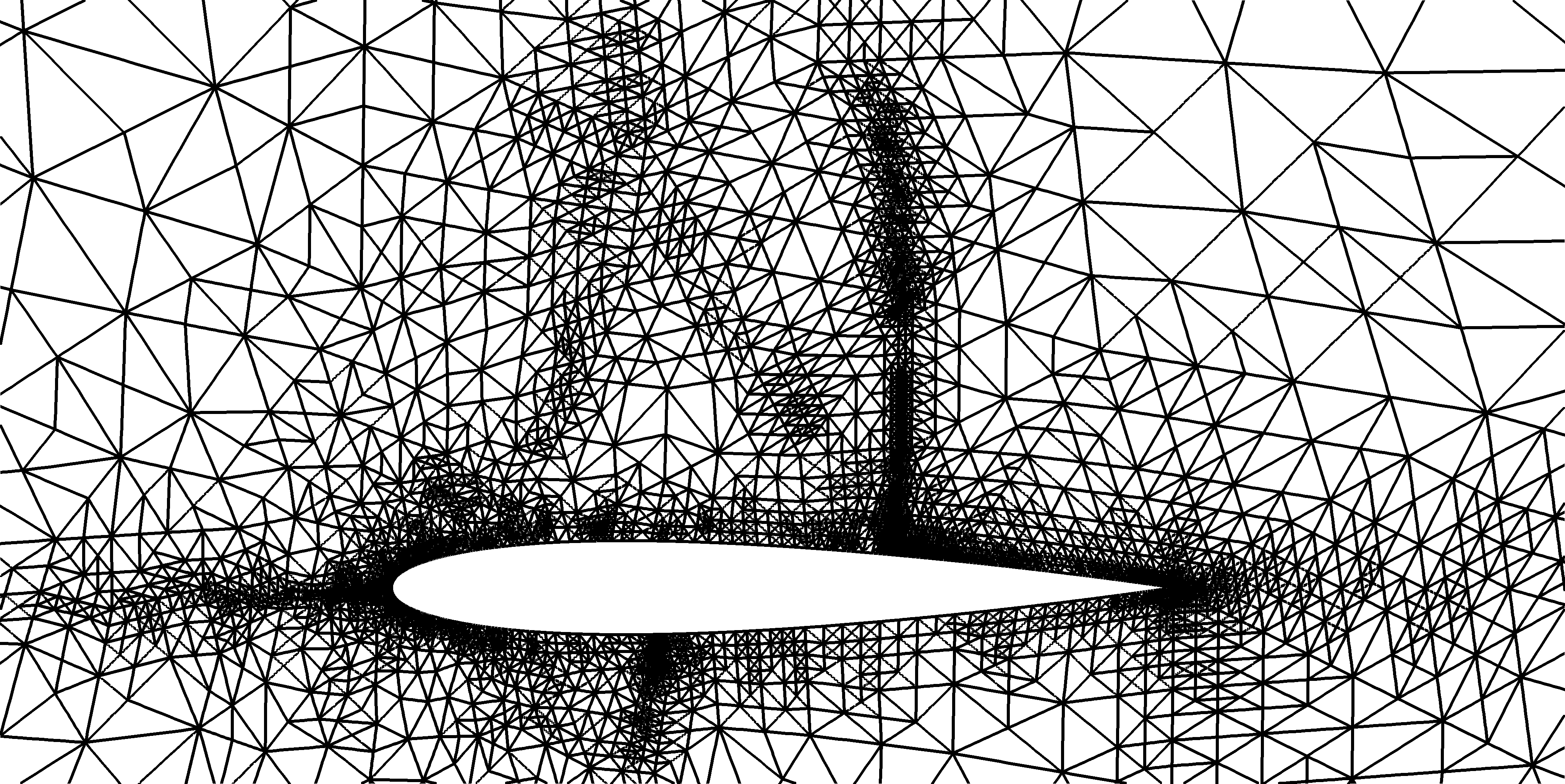}}
\caption{Meshes for the transonic Euler test case (${\rm Ma}_{\infty}=0.8$, $\alpha=1.25^{\circ}$)}
\end{figure}

The adjoint-based indicator performs well for all $p$. Even for this non-smooth flow (with a non-smooth dual solution) we approach the optimal rate of $2p+1$ (see Fig.~\ref{fig:naca_ma08_al125_ndof}). The residual-based adaptation is comparable during the first steps; then, however, the convergence gets worse and stagnates in some cases. Again, we can find the explanation for this by considering the refined meshes. The residual-based indicator focuses mainly on the upper shock (see Fig.~\ref{fig:naca_ma08_al125_mesh_h_res}). Both leading and trailing edge exhibit some refinement. The adjoint-based adapted mesh (see Fig.~\ref{fig:naca_ma08_al125_mesh_h}) shows a more complex refinement pattern. All of the regions mentioned at the beginning (leading and trailing edge, upper and lower shock) get refined. Thus, by properly weighting the local residuals with the adjoint, all important features become visible for the refinement process. Furthermore, the so-called lambda-feature (upstream of the primal shock) and the discontinuity emerging from the leading edge in the adjoint solution exhibit refinement (see Fig.~\ref{fig:naca_ma08_al125_adjoint2}). This serves for a proper representation of the adjoint solution which is necessary for both local and global error estimation.

The correction via the adjoint-based error estimate is again of minor use (see Fig.~\ref{fig:naca_ma08_al125_ndof}). This might be due to the complex features in both the primal and adjoint solution which are not resolved enough (see Fig.~\ref{fig:naca_ma08_al125_adjoint2}).

\FloatBarrier
\subsubsection{Laminar, Subsonic Flow}

We consider viscous subsonic flow around the NACA0012 airfoil with a free-stream Mach number of $\rm Ma_{\infty}=0.5$, a Reynolds number of $\rm Re=5000$ and an angle of attack $\alpha=1^{\circ}$. This flow is dominated by the relatively thin, laminar boundary layer.

\begin{figure}[h]
\centering
\subfloat[Mach number ($p=2$)\label{fig:naca_ma05_al1_re5000_ma}]{
\includegraphics[width=0.8\textwidth]{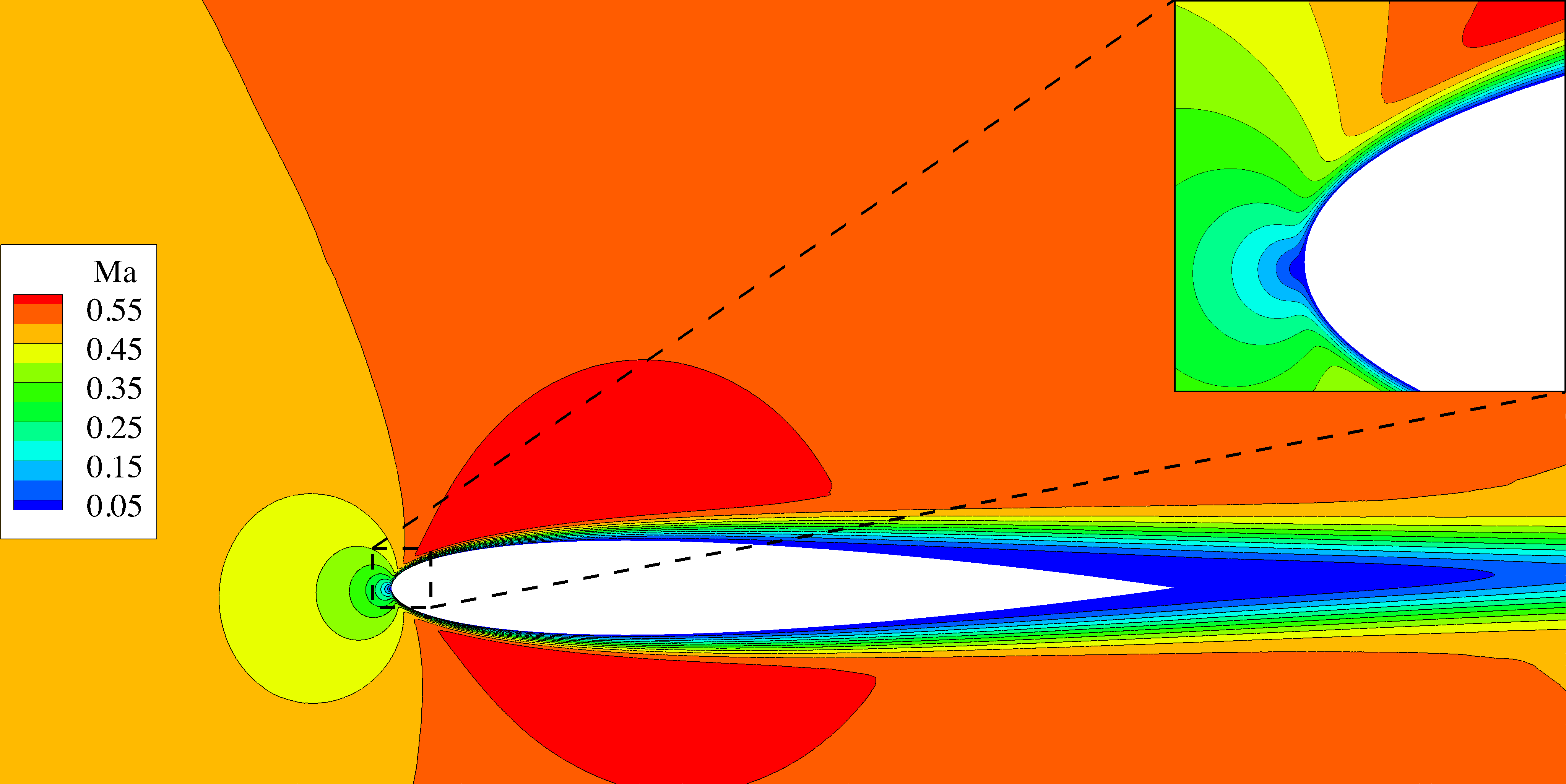}}\\
\subfloat[Adjoint x-momentum ($p=3$)\label{fig:naca_ma05_al1_re5000_adjoint2}]{
\includegraphics[width=0.8\textwidth]{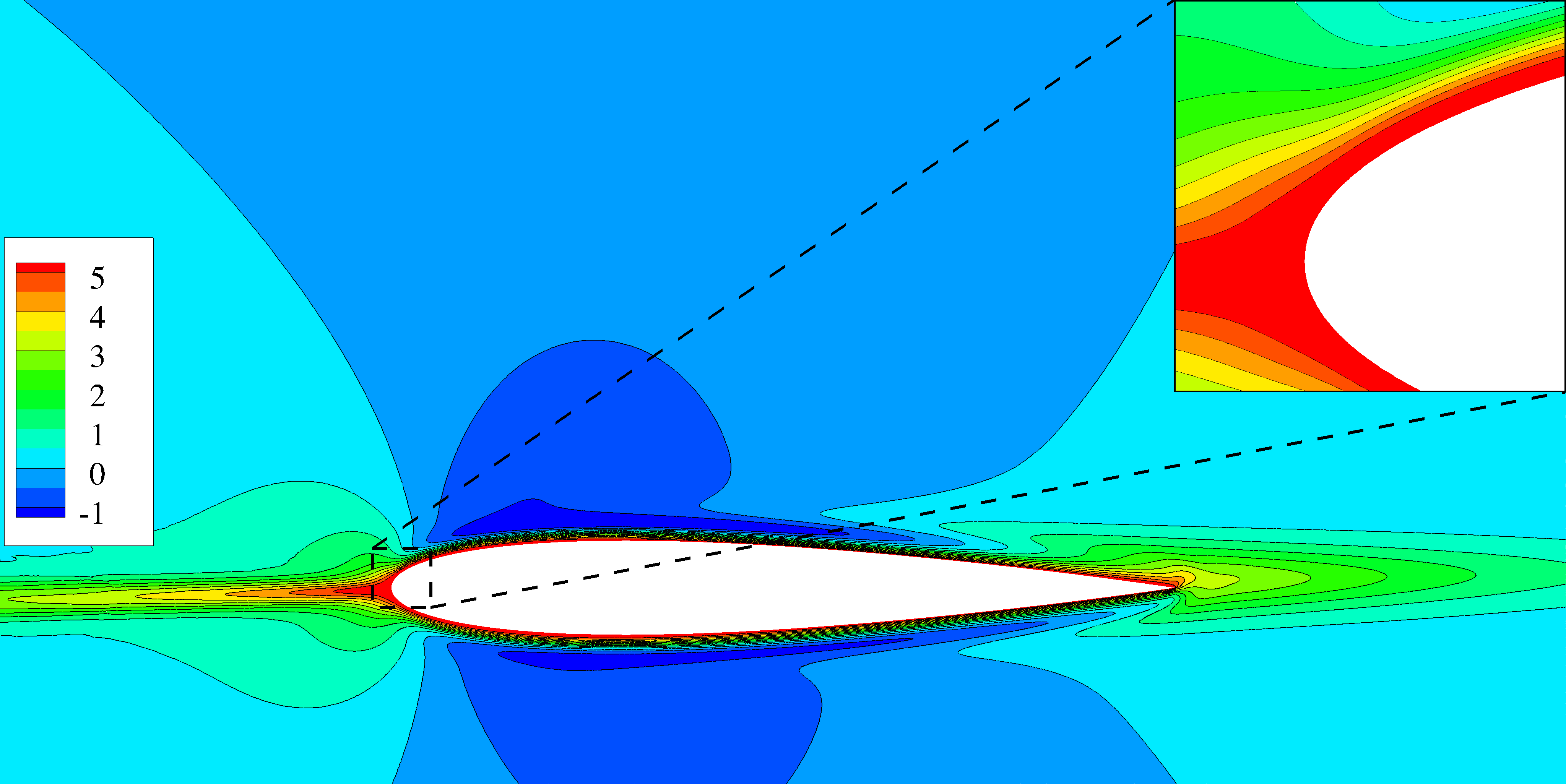}}
\caption{Primal and dual (drag) solution for the laminar test case (${\rm Ma}_{\infty}=0.5$, $\alpha=1^{\circ}$, $\rm Re=5000$).}
\end{figure}

In Fig.~\ref{fig:naca_ma05_al1_re5000_ma} and \ref{fig:naca_ma05_al1_re5000_adjoint2} contours of the Mach number and the adjoint x-momentum with respect to drag are given. The approximate adjoint solution is smooth due to the adjoint-consistent treatment of boundary conditions and target functionals. One can see that there is also a boundary layer evolving in the adjoint solution. The highest values in the adjoint are attained along the surface of the airfoil and upstream and downstream of leading and trailing edge, respectively.

\begin{figure}[h]
\centering
  \subfloat[$p=1$]{\includegraphics{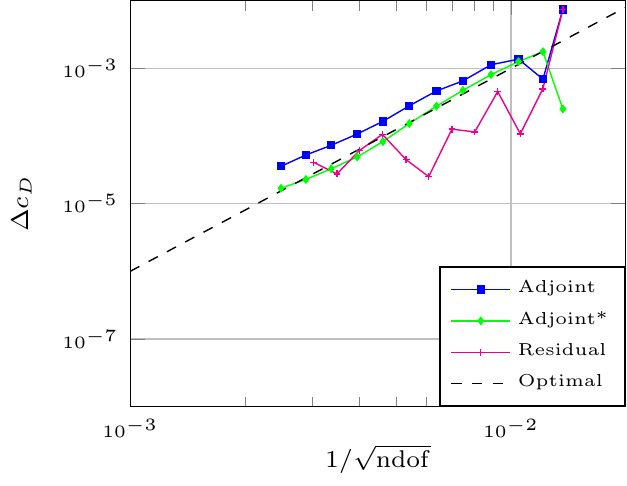}}\qquad
  \subfloat[$p=2$]{\includegraphics{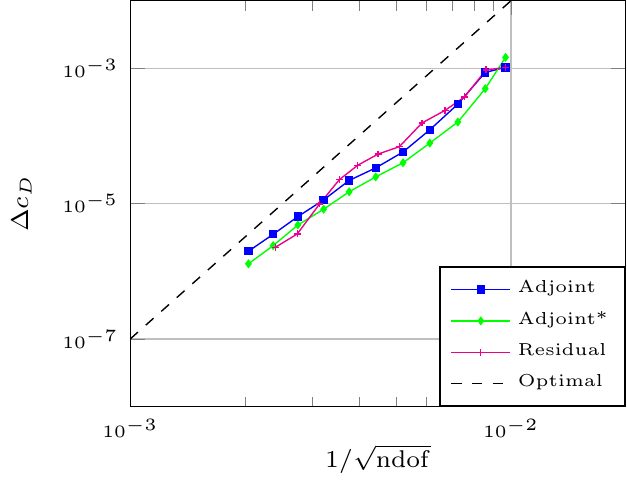}}\\
  \subfloat[$p=3$]{\includegraphics{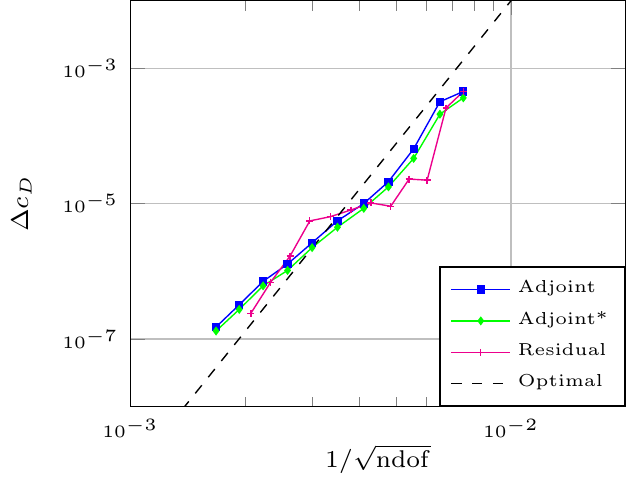}}\qquad
  \subfloat[$p=4$]{\includegraphics{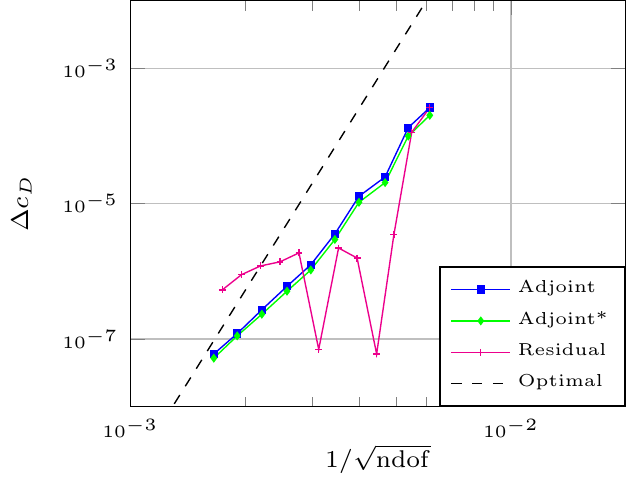}}
  \caption{Drag convergence with respect to degrees of freedom (${\rm Ma}_{\infty}=0.5$, $\alpha=1^{\circ}$, $\rm Re=5000$). A * denotes adjoint-based corrected values. Optimally, the error should converge as $\mathcal{O}\left(h^{2p+1}\right)$.}
  \label{fig:naca_ma05_al1_re5000_ndof}
\end{figure}

\begin{figure}[h]
\centering
\subfloat[Baseline mesh\label{fig:naca_ma05_al1_re5000_mesh_coarse}]{
\includegraphics[width=0.8\textwidth]{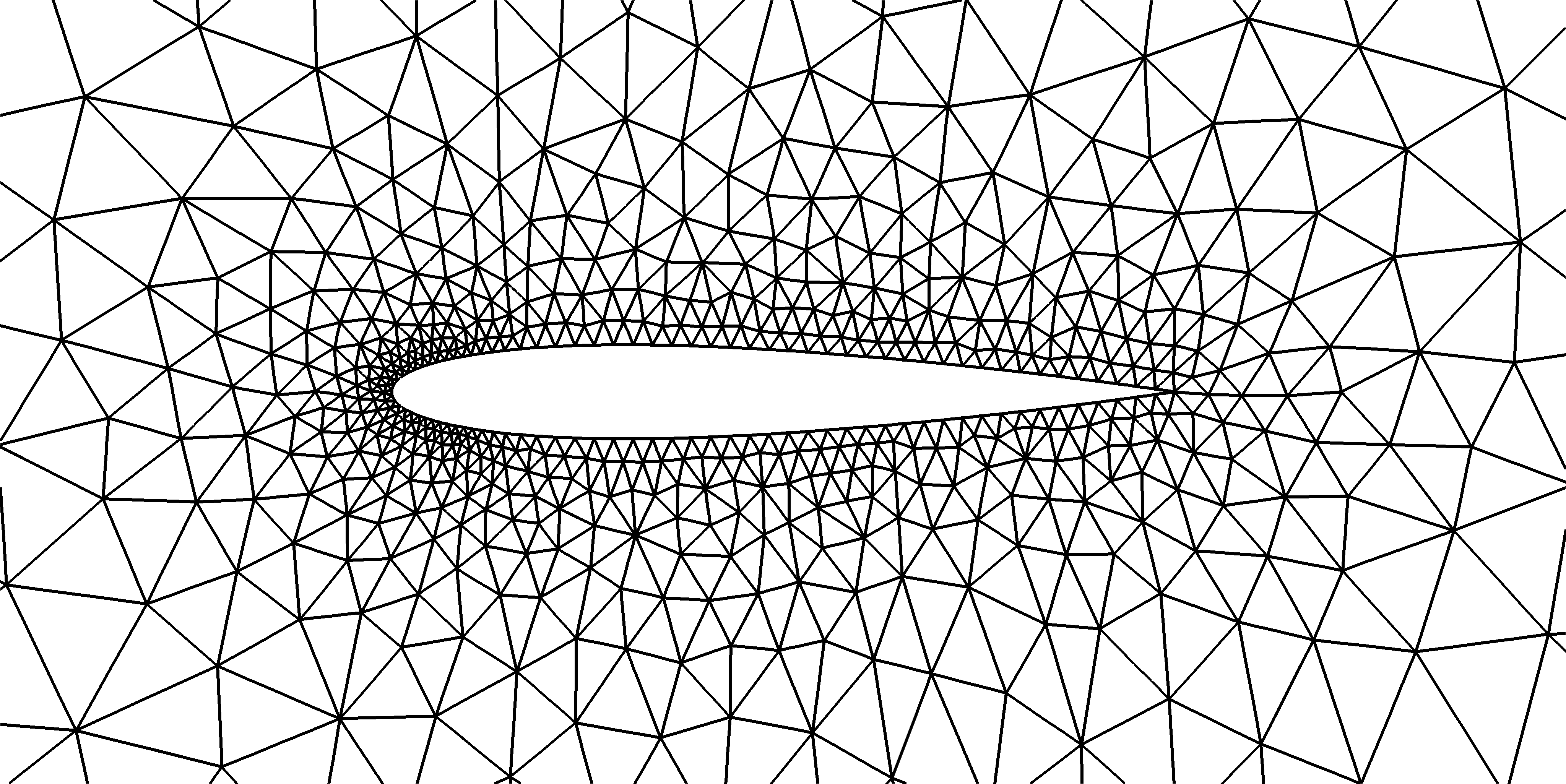}}\\
\subfloat[Residual-based refined mesh ($p=2$, $n_e\approx 18000$)
\label{fig:naca_ma05_al1_re5000_mesh_h_res}]{
\includegraphics[width=0.8\textwidth]{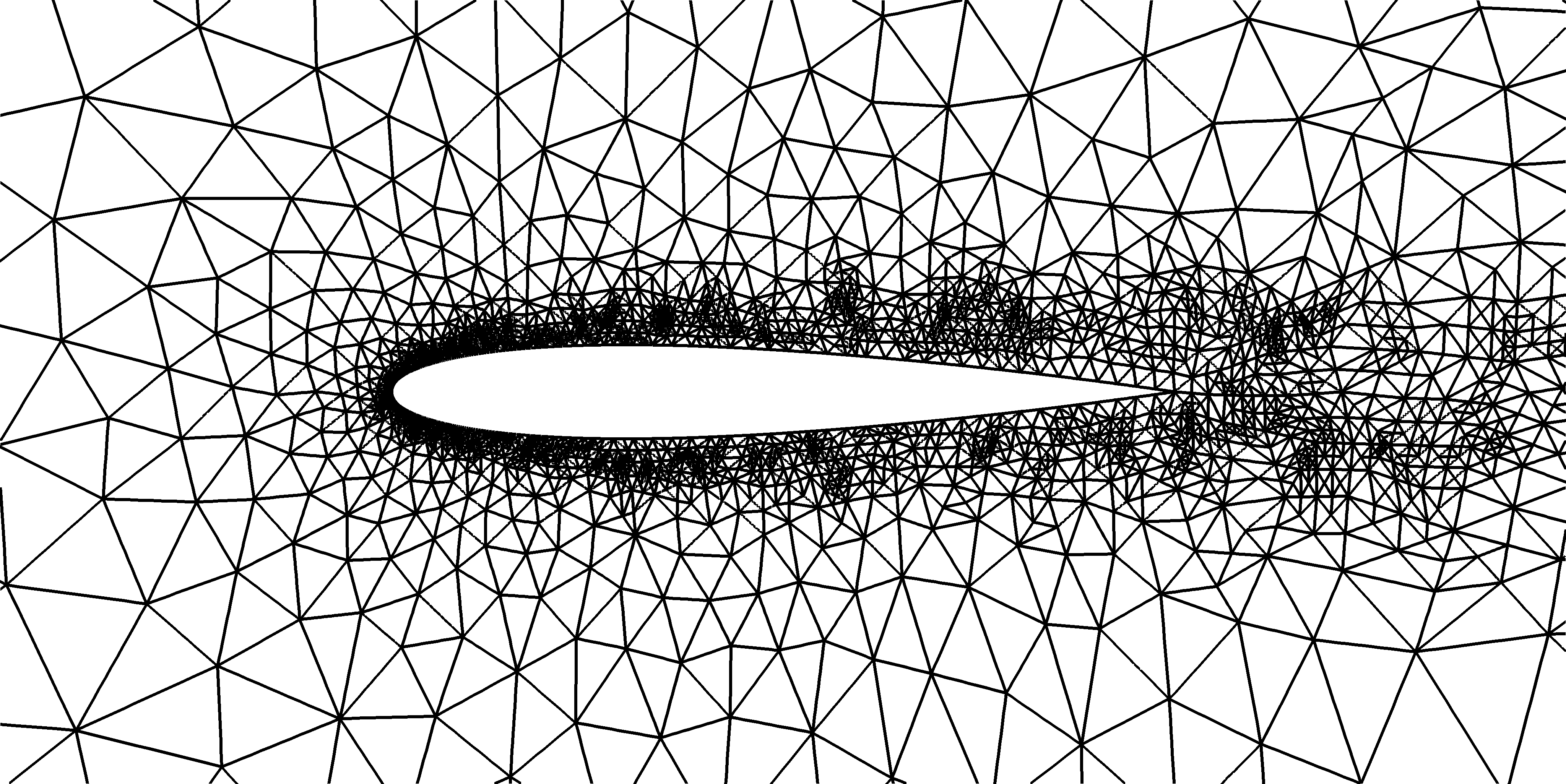}}\\
\subfloat[Adjoint-based refined mesh ($p=2$, $n_e\approx 18000$)
\label{fig:naca_ma05_al1_re5000_mesh_h}]{
\includegraphics[width=0.8\textwidth]{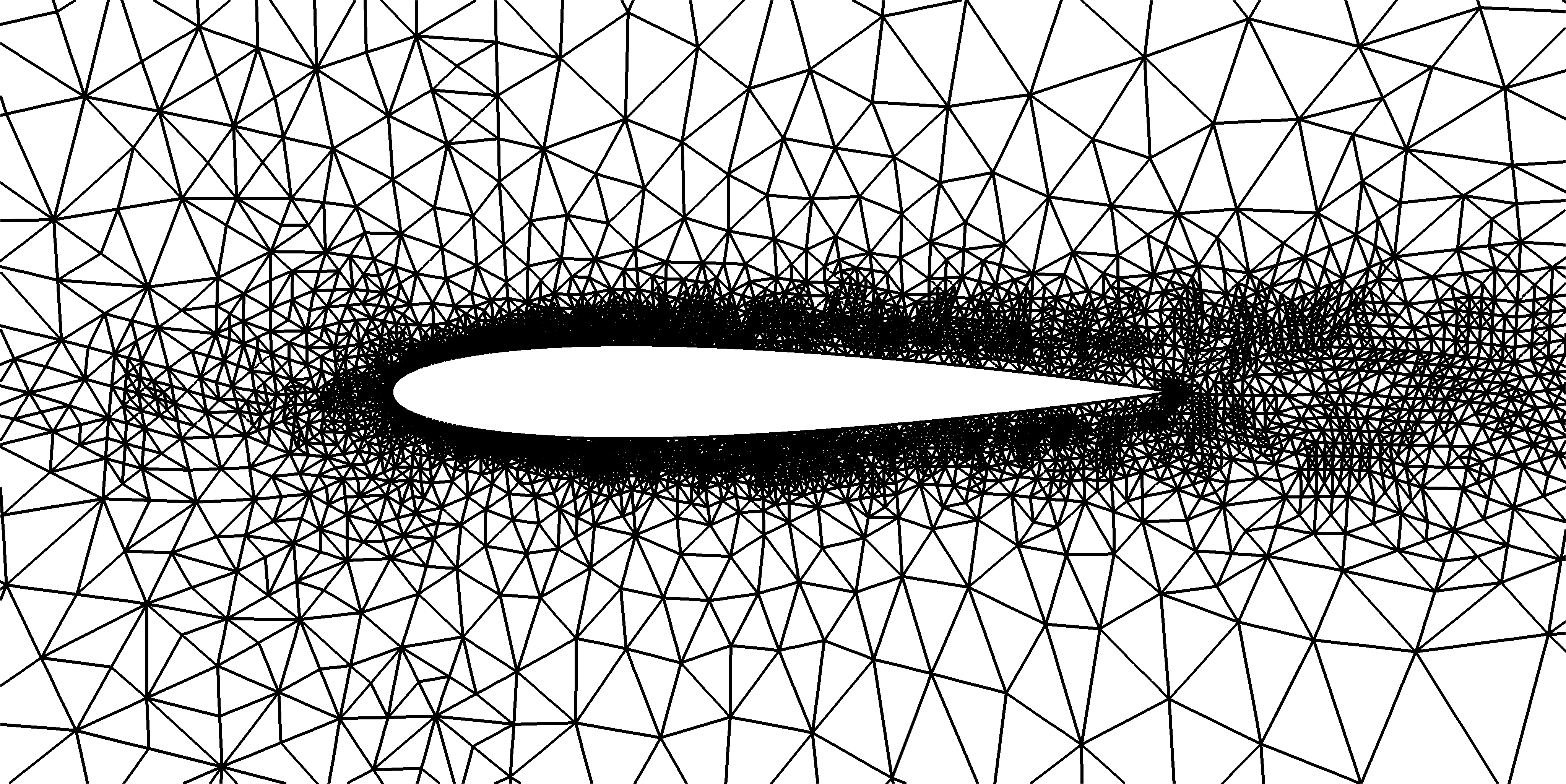}}
\caption{Meshes for the laminar test case (${\rm Ma}_{\infty}=0.5$, $\alpha=1^{\circ}$, $\rm Re=5000$)}
\end{figure}

$\theta=0.1$ showed good performance for the adaptive runs. In order to compute the error in the drag coefficient, a reference value was obtained on an adjoint-based adapted mesh with approximately $2.5\cdot 10^5$ degrees of freedom. Both adjoint- and residual-based adaptation perform comparably well (see Fig.~\ref{fig:naca_ma05_al1_re5000_ndof}). This is also visible in the refined meshes (see Fig.~\ref{fig:naca_ma05_al1_re5000_mesh_h_res} and Fig.~\ref{fig:naca_ma05_al1_re5000_mesh_h}). They both show large refinements in the boundary layer. Only for $p=1$ and $p=4$ the convergence path of the residual-based refinement shows some oscillatory behavior. The adjoint-based indicator yields a reduction of the error in every step. The optimal rate of $2p+1$ can however not be obtained. This might be due to our type of refinement. For test cases which exhibit strongly anisotropic features (like the boundary layer in this one), anisotropic mesh adaptation should be used to yield an appropriate mesh.

The correction via the adjoint-based error estimate results again only in a slight reduction of the error (see Fig.~\ref{fig:naca_ma05_al1_re5000_ndof}). Like in the transonic test case, this might be due to the still unresolved features in the adjoint solution (see Fig.~\ref{fig:naca_ma05_al1_re5000_adjoint2}).

\section{Conclusion and Outlook}
\label{sec:conclusions}
We presented an adjoint-based mesh adaptation methodology for hybridized discontinuous Galerkin methods. The scheme  not only decreases the number of globally coupled degrees of freedom by hybridization, but also distributes them in an efficient way. The adjoint-based adaptation strategy proved to be superior to both uniform and residual-based mesh refinement. In nearly all cases, super-convergence in target functionals could be recovered with the aid of adjoint-based adaptive refinement. In the case of anisotropic solution features, isotropic mesh refinement was found to be suboptimal. Here, an adequate anisotropic mesh adaptation procedure should be employed. This is planned for future work. Furthermore, the convergence behavior of the global error estimate should be investigated in more detail.

In a forthcoming paper we will extend our method to hp-adaptation. (see \cite{balan2013hybridadjointhp} for preliminary results.) In a concurrent effort, we have compared our method with a standard DG method in terms of efficiency \cite{woopen2013comparison}. Finally, we plan to extend our computational framework to three dimensional problems \cite{woopen2014hdg3d}. Then, adaptivity will play an even more crucial role, as the problem size increases drastically compared to the two dimensional case.

\section*{Acknowledgments}

Financial support from the Deutsche Forschungsgemeinschaft (German Research Association) through grant GSC 111, and by the Air Force Office of Scientific Research, Air Force Materiel Command, USAF, under grant number FA8655-08-1-3060, is gratefully acknowledged.

\bibliographystyle{plainnat}

\end{document}